\documentclass{amsart}
\usepackage{times}

\usepackage{amsfonts}
\usepackage{bbold}
\usepackage{amscd}
\usepackage{amssymb}
\usepackage{latexsym}
\usepackage{amsthm}
\usepackage{amsxtra}

\usepackage[bookmarks,colorlinks,pagebackref]{hyperref} 
\input{xy}
\xyoption{all}

\theoremstyle{plain}
\newtheorem{stelling}{Theorem}

\renewcommand\thestelling{\textnormal{\Alph{stelling}}}

\swapnumbers
\newtheorem{theorem}[subsection]{Theorem}
\newtheorem{corollary}[subsection]{Corollary}
\newtheorem{lemma}[subsection]{Lemma}
\newtheorem{proposition}[subsection]{Proposition}

\theoremstyle{definition}
\newtheorem{definition}[subsection]{Definition}

\theoremstyle{remark}

\swapnumbers

\newcommand{\emptyprop}{q}

\newcommand \acf{algebraically closed field}

\newcommand \binomial[2]{{\bigl( \begin{matrix} #1\cr#2\cr\end{matrix} \bigr)}}
\newcommand \ch{characteristic}
\newcommand \CM{Coh\-en-Mac\-au\-lay}
\newcommand \complet[1]{\widehat {#1}} 
\newcommand \DVR{discrete valuation ring}

\newcommand \homo{homomorphism}

\renewcommand\iff{if and only if}

\newcommand \iso{\cong}
\newcommand \loc{{\mathcal {O}}}
\newcommand \los{\L os' Theorem}
\newcommand \map[1]{{\newcommand{\tmpprop}{#1q}  \if\tmpprop\emptyprop \to\else \xrightarrow{{\phantom{i}{#1}\phantom{i}}}\fi}} 
\newcommand \maxim{\mathfrak m}
\newcommand \nat{\mathbb N}
\newcommand \norm[1]{\left|#1\right|}

\newcommand \op\operatorname
\newcommand \pol[2]{#1[#2]}
\newcommand \pow[2]{#1[[#2]]}

\newcommand \pr{\mathfrak p}
\newcommand \primary{\mathfrak g}

\newcommand \restrict [2]{\left.#1\right|_{{#2}}}
\newcommand \rij[2]{(#1_1,\dots,#1_{#2})}
\let\sub\subseteq

\newcommand \tensor{\otimes}
\newcommand \tor[4]{\operatorname{Tor}^{#1}_{#2}(#3,#4)}

\newcommand \zet{\mathbb Z}

\newcommand{\commdiagram}[9][]{%
\begin{equation}
{\newcommand{\tmpprop}{#1q} 
\if\tmpprop\emptyprop \relax\else \label{#1}\fi}
\begin{aligned}%
\mbox{
\begin{picture}(130,90)%
\put(120,70){\vector( 0,-1){50}}%
\put(10,80){\vector( 1, 0){100}}%
\put(0,70){\vector( 0,-1){50}}%
\put(10,10){\vector( 1, 0){100}}%
\put(115,80){\makebox(0,0)[l]{$#4$}}%
\put(5,80){\makebox(0,0)[r]{$#2$}}%
\put(115,10){\makebox(0,0)[l]{$#9$}}%
\put(5,10){\makebox(0,0)[r]{$#7$}}%
\put(-3,50){\makebox(0,0)[r]{$#5$}}
\put(123,50){\makebox(0,0)[l]{$#6$}}
\put(60,3){\makebox(0,0)[c]{$#8$}}
\put(60,88){\makebox(0,0)[c]{$#3$}}
\end{picture}}
\end{aligned}
\end{equation}}

\newcommand\commtrianglefront[7][]{%
\begin{equation}
{\newcommand{\tmpprop}{#1q} 
\if\tmpprop\emptyprop \relax\else \label{#1}\fi}
\begin{aligned}%
\mbox{
\begin{picture}(120,80)%
\put(55,68){\vector(-1,-2){30}}
\put(65,68){\vector(1,-2){30}}
\put(30,5){\vector(1,0){60}}
\put(60,75){\makebox(0,0)[c]{$#2$}}
\put(25,5){\makebox(0,0)[r]{$#4$}}
\put(95,5){\makebox(0,0)[l]{$#6$}}
\put(60,0){\makebox(0,0)[c]{$#5$}}
\put(37,43){\makebox(0,0)[r]{$#3$}}
\put(83,43){\makebox(0,0)[l]{$#7$}}
\end{picture}}
\end{aligned}
\end{equation}}

\newcommand\commtriangleback[7][]{%
\begin{equation}
{\newcommand{\tmpprop}{#1q}
\if\tmpprop\emptyprop \relax\else \label{#1}\fi}
\begin{aligned}%
\mbox{
\begin{picture}(120,80)%
\put(55,70){\vector(-1,-2){30}}
\put(65,70){\vector(1,-2){30}}
\put(30,5){\vector(1,0){60}}
\put(60,75){\makebox(0,0)[c]{$#2$}}
\put(25,5){\makebox(0,0)[r]{$#6$}}
\put(95,5){\makebox(0,0)[l]{$#4$}}
\put(60,0){\makebox(0,0)[c]{$#5$}}
\put(37,43){\makebox(0,0)[r]{$#7$}}
\put(83,43){\makebox(0,0)[l]{$#3$}}
\end{picture}}
\end{aligned}
\end{equation}}

\newcommand \ul[1]{#1\mathstrut_\natural}  
\newcommand \ulsep[1]{#1\mathstrut_\sharp} 
\newcommand \tuple[1]{\mathbf{#1}} 
\newcommand \defpair{deformator} 
\newcommand \tuplepow[2]{\tuple{#1}^{(#2)}} 
\newcommand \poinc[2]{\op{Poin}^{\op{#2}}(#1)} 

\newcommand \seq[2]{#1\mathstrut_{#2}}
\newcommand \app[2]{J^{#2}#1}

\newcommand \scal[2]{\complet{#1}^{#2}}

\newcommand\near[1]{\ulcorner #1\lrcorner}

\newcommand  \deform{jet}
\newcommand  \Deform{Jet}
\newcommand  \infdef{\deform}
\newcommand  \simloc{\mathbb{Sim}}
\newcommand  \simlocdef{\mathbb{SimDef}}
\newcommand  \simlocpar{\mathbb{SimPar}}
\newcommand  \simcm {\mathbb{SimCM}}
\newcommand  \siman {\mathbb{SimAn}}
\newcommand \cau[1]{\op{Cau}(#1)}
\newcommand  \ball{\mathbb B}

\renewcommand \inf[1]{\op{Inf}(#1)}

\renewcommand  \c{ca\-ta}

\newcommand \ac[1]{#1^{\text{alg}}}

\title {Classifying singularities up to analytic extensions of scalars}
\author{Hans Schoutens}
\date{\today}
\thanks{Partially supported by the National Science Foundation and a PSC CUNY
grant}

\begin{document}

\begin{abstract} 
The singularity space consists of all germs $(X,x)$, with $X$
a Noetherian scheme and $x$ a point, where we identify two such germs if they
become the same after an analytic extension of scalars.
 This is a Polish space for the metric given by the order
to which infinitesimal neighborhoods, or jets, agree after
base change.  In other words, the classification of singularities up to
analytic extensions of scalars is a smooth problem in the sense of descriptive
set-theory.  Over $\mathbb C$, the following two classification problems up to isomorphism are now smooth:  (i) analytic germs; and (ii)   polarized schemes.\end{abstract}

\maketitle

\section{Introduction}
Roughly speaking, a classification problem consists of a class of objects
together with an equivalence relation telling us which
objects to identify; a solution to this problem is then an
`effective' or `concrete' description of the quotient, preferably by a
`system of complete invariants'. What constitutes a reasonably concrete or
effective solution to a classification problem, however,  might depend on one's
purposes or even one's taste. Descriptive set-theory proposes
smoothness to be the decisive indication that a classification is explicit and/or
concrete (see for instance \cite{HKL,Hjo} for a discussion). More precisely,
recall that a \emph{Polish} space is  a complete metric space  containing a
countable dense subset. Considering a Polish space to be concrete is justified by 
the fact that its underlying Borel structure   is
in essence equal to the  standard  Borel space $\mathbb R$. With this in mind, an
equivalence relation on a Polish space, and by extension, the classification problem it encodes, is called \emph{smooth} if there is a Borel map to a Polish space which factors through the quotient. A more suggestive, albeit slightly less precise formulation is that, up to a Borel isomorphism, equivalence classes are completely classified by real numbers.  

Most  classification problems in algebraic geometry, like
classifying varieties over a fixed \acf\  up to isomorphism or up to bi-rational
equivalence, are
not known to be smooth. Of course, this is in no way preventing geometers to
seriously, and often successfully, work on   these classification problems.
It
would be nice to know though what their descriptive set-theoretical status is.
 In this paper, I will propose a \emph{local} classification problem, which
will fall at the right side of the dividing line:    one can `concretely', that is to say, smoothly, classify
germs of points on arbitrary Noetherian schemes up to similarity (a slightly weaker equivalence relation than
 the isomorphism relation). Using this general result, we can also deduce some  smoothness results for certain isomorphism problems. For  \emph{analytic germs}, that is to say,   formal completions of  germs (in the sense of \cite[II.9]{Hart}), we have: 
  
 \begin{theorem}\label{T:formgerm}
The classification,  up to isomorphism, of analytic germs over   an \acf\ of size the continuum, is smooth.
\end{theorem}

This also enables us to obtain  a  smooth classification problem of a more global nature, namely for projective schemes together with a choice of a very ample line bundle, the so-called \emph{polarized schemes}.

\begin{theorem}\label{T:pol}
The classification,  up to isomorphism, of polarized schemes  over an \acf\ of size the continuum,   is smooth.
\end{theorem}
 
For the proof of our main smoothness result, we associate to a point its local ring, thus reducing the problem   to the study
of the category of all Noetherian local rings. If we were to classify
these only up to isomorphism, then 
as part of this problem, we would have to
classify already all fields, and even for countable fields \cite{FriSta} or
fields of finite transcendence degree \cite{ThoVel} these are non-smooth
problems. Hence to circumvent this arithmetical obstruction, we can either
fix the residue field---the route taken for the two isomorphism problems stated above---or, otherwise,  allow
for `extensions of scalars', resulting in the 
identification of any two fields of the same \ch.   Even after taking the latter modification, the local  classification problem is   probably 
still not smooth. We   introduce one further identification,  inspired
 in part by Grothendieck's suggestion that one should consider working in  the 
etale topos instead of the  (classical) Zariski topos. A down-to-earth interpretation of
this point of view is that two local rings can be considered identical if they
have a common etale
extension, or more generally, if they have the same completion. In summary, we 
say that two Noetherian local rings are \emph{similar} if they can be made isomorphic by
an \emph{analytic extension of scalars}, that is to say, by the process of
extending scalars and taking completion.
To also make sense of this in mixed \ch, we subsume these types of extensions 
under the larger class consisting of all formally etale (=unramified and faithfully
flat) extensions.
We will show that \emph{similar} points (meaning that their corresponding
local rings are similar) have the same type of
singularity (see Theorem~\ref{T:transsim}). As a spinoff of this investigation,
we obtain a flatness criterion generalizing a result
of Koll\'ar:

\renewcommand\thestelling{{\ref{T:ff}}}
\begin{stelling}
 Let  $R\to S$ be a local \homo\ between Noetherian
local rings and suppose $R$ is  an excellent normal   domain with perfect residue
field. If  $\op{dim}(R)=\op{dim}(S)$ and $R\to S$ is unramified, then $R\to S$ is
faithfully flat.
\end{stelling}

Our assertion that classifying points up to similarity is smooth is 
established by effectively putting a metric on the space of similarity classes
$\simloc$, called the \emph{\deform\ metric}. We will prove that the induced
topology is complete, and that the collection of similarity classes of Artinian
local rings with a finitely
generated residue field is a countable dense subset. This shows that $\simloc$ is
a Polish space and hence classification up to similarity is a smooth problem.
The \deform\ metric on $\simloc$ is induced by a semi-metric on the class of all 
Noetherian local rings. In terms of (germs of) points, this semi-metric measures
to which order the \deform{s}  of two points agree.  In fact, the proof yields that  for   classification up to
similarity, the
collection of all  \deform{s} of a point form a
complete set of invariants.

So far, all concepts are algebraic-geometric in
nature, but the existence of limits relies on a tool from model-theory, to wit,
the
ultraproduct construction. Of course, the ultraproduct of Noetherian local rings
is in general no longer
Noetherian. However, if we have a Cauchy sequence of Noetherian
local rings, then their \emph{\c{}product}, obtained by killing
  all infinitesimals in the ultraproduct, yields a
complete Noetherian local ring, which, up to similarity, is the limit of the
sequence. 
 
%
%

\section{Limits and ultraproducts}

Let $(\Sigma,d)$ be a \emph{semi-metric space}. In this paper, we understand
this to mean that the semi-metric is \emph{non-archimedean}, 
 that is to say,   $d(x,y)\leq\max\{d(x,z),d(y,z)\}$ for all $x,y,z\in\Sigma$, 
and \emph{bounded}, that is to say, after possibly normalizing the metric, $d(x,y)\leq 1$ for all $x,y\in\Sigma$.
We call $d$ a
\emph{metric}, if $d(x,y)=0$ \iff\ $x=y$. To include the \deform\
metric in our
treatment, we allow for $\Sigma $ to be merely a class. We say that two elements
$r,s\in\Sigma$ are   \emph{$d$-equivalent},
written $r\sim_d s$, if $d(r, s)=0$.   The quotient space  $\Sigma/\sim_d$
has an induced semi-metric which is in fact a metric; we therefore call this
quotient the   \emph{metrization} of $(\Sigma,d)$.

Let $(\seq \Sigma w,\seq dw)$  be semi-metric spaces, for $w\in\nat$. We will
identify the elements of the   product $\Pi:=\prod_w\seq \Sigma
w$ with the sequences $\mathbf r\colon\nat\to \Pi$ such that $\mathbf r(w)\in\seq\Sigma
w$ for each $w$. The \emph{product semi-metric} on $\Pi$
is given by letting the distance $d(\mathbf r,\mathbf s)$ between two  
sequences  
$\mathbf r$ and $\mathbf s$ be the lim-inf of the distances
$\seq dw(\mathbf r(w),\mathbf s(w))(\leq 1)$ of their respective components.  
Below, we will introduce weaker semi-metrics on $\Pi$, induced by ultrafilters.

\subsection*{Cauchy sequences}
Let $\mathbf r$ be a sequence  in $\Sigma$ (meaning that all $\mathbf
r(w)\in\Sigma$) and let $\mathbf r^+$ be its \emph{twist}, given as the
sequence whose $w$-th element is $\mathbf r(w+1)$. We call $\mathbf r$ a
\emph{Cauchy sequence} if $\mathbf r\sim\mathbf r^+$ (with respect to the product
semi-metric). One verifies that $\mathbf r$  is a 
Cauchy sequence, if for each $\varepsilon >0$, there exists an
$N$ such that  $d(\mathbf r(w),\mathbf r(v))< \varepsilon$ for all $v,w>N$, 
and that two Cauchy sequences  $\mathbf r$ and $\mathbf s$ are equivalent
 if for 
each $\varepsilon>0$, there exists an $N$ such that
$d(\mathbf r(w),\mathbf s(w))<\varepsilon$ for all $w>N$. Let $\cau{\Sigma,d}$, or
simply, $\cau \Sigma $, denote  the set of all Cauchy sequences in
$\Sigma $ with the induced product semi-metric. There is a natural isometry $\Sigma\to
\cau \Sigma $ sending $x$ to the constant sequence $\mathbf x$ given as
$\mathbf x(w):=x$; we will identify the element $x$ with its constant sequence in
$\cau \Sigma $.

A \emph{limit} of a sequence $\mathbf r$ is an element $x\in \Sigma $ such that $\mathbf r\sim x$.  It is easy to see that if
$\mathbf r$ has  a limit, then it must be Cauchy. We call $(\Sigma,d)$
\emph{complete} if every Cauchy sequence has a unique limit. This implies in
particular that  $d$ is a metric.  We define the \emph{completion} of
$(\Sigma,d)$ as the
metrization $\complet \Sigma:= \cau \Sigma/\sim$ of the semi-metric space $\cau\Sigma$;
it is a complete metric space containing $\Sigma $ as a dense subspace.

\subsection*{Adic metric}
A local ring $(R,\maxim)$ comes with a canonical  semi-metric, its
\emph{$\maxim$-adic semi-metric} defined as follows: the \emph{order} of an
element
$x\in R$ is the  supremum of all $n$ for which $x\in\maxim^n$;  the distance $d_R(x,y)$
 between two elements is then equal to $2^{-n}$ where $n$ is the order of
$x-y$ (we allow the case $n=\infty$, with the convention that $2^{-\infty}=0$).
The subset of all elements which are $d_R$-equivalent to zero forms 
an ideal, equal to the intersection of all the powers $\maxim^n$; it is called
the
\emph{ideal of infinitesimals} of $R$ and is denoted  
$\inf R$. By Krull's intersection theorem, if $R$ is
Noetherian, then $\inf R=0$.
The completion
of $R$ in the $\maxim$-adic semi-metric will be denoted $\complet R$. If $R$ has finite
embedding dimension, 
 then $\complet R$ is a
complete Noetherian local
ring by \cite[Theorem 2.2]{SchFinEmb}. 

Below, we will define a semi-metric on the class of all Noetherian local
rings, not to be confused with the adic metric on a single Noetherian local
ring. To calculate limits in the former semi-metric, we need   a notion from
model-theory:
the ultraproduct construction (some references for ultraproducts are    
\cite{EkUP,Hod,Roth,SchUlBook}, or the brief review in \cite[\S2]{SchNSTC}) .

\subsection*{Ultraproducts and \c{}products}
Let $(\seq Rw,\seq\maxim w)$, for $w\in\nat$, be a sequence  of
Noetherian local rings. Let $\mathcal U$ be an ultrafilter on
$\nat$, which we always assume to be non-principal. The
\emph{ultraproduct} of  the $\seq Rw$ with respect to $\mathcal U$, denoted
$\ul R$, is   obtained from the product $\Pi:=\prod_w\seq Rw$  by modding out the ideal of all sequences almost all of whose
entries are zero (it is customary to use the expression ``almost all'' to mean
``all indices belonging to a member of the ultrafilter''). The particular choice\footnote{\label{f}There is really no reason to restrict
only to ultraproducts on a countable index set, although it is the only type we
will use in this paper. However, for the
\c{}product (see below) to be Noetherian and complete, we do have to impose that
the ultrafilter is countably incomplete, which automatically holds on
countable index sets and always exists on arbitrary index sets.} of ultrafilter
$\mathcal U$ does not matter for our purposes, and hence we do not include it in our notation. Although not useful for proving results, let me recall an alternative construction from \cite[Theorem 2.5.4]{SchUlBook}: there exists a minimal prime ideal $\primary$ of the Cartesian power $\zet^\nat$, containing the direct sum ideal $\zet^{(\nat)}$, such that $\ul R=\Pi/\primary\Pi$, where we view the Cartesian product $\Pi$ as an algebra over $\zet^\nat$ in the natural way; and conversely, any such prime ideal determines in this way an ultraproduct of the $\seq Rw$.

The ultraproduct $\ul R$ is again a
local ring, with maximal ideal $\ul\maxim$  given as
the ultraproduct of the $\seq\maxim
w$. In general, however, $\ul R$ will no longer be
Noetherian. If almost all $\seq
Rw$ have embedding dimension at most $n$, then so does $\ul R$. A key role will played by the homomorphic image of $\ul R$ modulo its ideal of infinitesimals $\inf{\ul R}$, which we call the \emph{\c{}product} of the $\seq
Rw$ and which we denote by $\ulsep R$.  
A more direct way 
for defining the \c{}product, although less useful in proofs, is as
follows: on the product $\Pi$, the ultrafilter $\mathcal U$ induces
a semi-metric $d_{\mathcal U}$ by the condition  that $d_{\mathcal U}(\mathbf r,\mathbf
s)\leq \varepsilon$ for some $\varepsilon$ \iff\  $d_{\seq Rw}(\mathbf r(w),\mathbf s(w))\leq
\varepsilon$ for almost all $w$. The
\c{}product is then the metrization of $(\Pi, d_{\mathcal U})$,
that is to say, the residue ring of the product modulo the ideal of all sequences
which are $d_{\mathcal U}$-equivalent to zero.

If almost all $\seq Rw$
have
embedding dimension at most $n$, then so does the cataproduct $\ulsep R$.
Moreover, by the saturatedness  property of ultraproducts, the \c{}product is 
a complete local ring, whence Noetherian by
\cite[Theorem 29.4]{Mats} (for more details see  \cite[Lemma 5.6]{SchFinEmb} or \cite[Theorem 12.1.4]{SchUlBook}). The same argument also shows that the  $\seq Rw$ and their completions 
$\seq{\complet R}w$ have the same  
\c{}product.

We will only consider \c{}products of Noetherian
local rings of bounded embedding dimension, so that we tacitly may assume that
they are complete and Noetherian. 
In case all $\seq Rw$ are equal to a fixed 
Noetherian local ring $R$, then their  ultraproduct $\ul R$ and \c{}product $\ulsep R$ are called, respectively, the  \emph{ultrapower} and \emph{\c{}power}
of $R$. By \los, ultrapowers commute with base change, that is to say
$\ul{(R/I)}\iso \ul R/I\ul R$; the same is true for \c{}powers  by \cite[Corollary
5.7]{SchFinEmb}:

\begin{lemma}\label{L:ulsepbc}
If $R$ is a Noetherian local ring and $I$ an ideal in $R$, then $\ulsep{ (R/I)}=\ulsep
R/I\ulsep R$.
\end{lemma}

\section{Scalar extensions}\label{s:scext}

Cohen's structure theorems for complete Noetherian local rings will play an
essential
role in this paper, so we quickly review the relevant properties; a good
reference for all this is \cite[\S29]{Mats}.

\subsection*{Cohen's structure theorems} 
For each field $\kappa$ of prime \ch\ $p$, there exists a unique complete \DVR\
$V$ of \ch\ zero 
whose residue field is $\kappa$ and whose maximal ideal is $pV$; we call $V$ the
\emph{complete $p$-ring} over $\kappa$.
Let $R$ be
a Noetherian local ring with residue field $\kappa$. We say that $R$ has \emph{equal
\ch} if $R$ and $\kappa$ have the same \ch; in the remaining case, we say that $R$
has  \emph{mixed \ch}. Assume $R$ is moreover complete and let $X$ be a finite
tuple of indeterminates.  Cohen's structure theorems now claim, among other
things, the following:
\begin{itemize}
\item if $R$ has equal \ch, then it is a homomorphic image of $\pow \kappa X$;
\item if $R$ has mixed \ch, then it is a homomorphic image of $\pow VX$,
where $V$ is the complete $p$-ring over $\kappa$.
\end{itemize}

Let $(R,\maxim)$ be a Noetherian local
with residue field $\kappa$ and let $\lambda$ be a field extension of $\kappa$. With a
\emph{scalar extension of $R$ over $\lambda$} we mean 
a local $R$-algebra 
$(S,\mathfrak n)$ with residue field $ \lambda $ such that $R\to S$ is faithfully flat, $\mathfrak
n=\maxim S$ and $R\to S$
induces the embedding $\kappa\subseteq \lambda$ on the residue fields. A \emph{scalar extension} of a local ring $R$ is
then a scalar extension of $R$ over some
field extension of its residue field. The condition that $\mathfrak n=\maxim S$
is also expressed by saying that $R\to S$  has \emph{trivial closed fiber} or that
it is \emph{unramified}.  In other words, a scalar extension is the same as an
unramified, faithfully flat \homo\ (also called a \emph{formally etale} extension). By
\cite[$0_{III}$ 10.3.1]{EGA}, for any Noetherian local ring $R$ and any extension $l$ of its residue field, at least one
scalar extension of $R$ over $l$ exists; we will reprove this in
Corollary~\ref{C:scal} below.

\begin{proposition}\label{P:scalfact} 
Consider the following commutative triangle of local
\homo{s}
between Noetherian local rings
\commtrianglefront {(R,\maxim)} {f} {(S,\mathfrak n)} {g}{(T,\pr)}  {h} 
If any two   are scalar extensions, then so is the third.
\end{proposition}
\begin{proof} 
It is clear that the composition of two scalar extensions is again
scalar. Assume $g$ and $h$ are scalar
extensions. Then $f$ is faithfully flat and $\maxim T=\pr=\mathfrak nT$. Since $g$ is faithfully flat, we  get
$\maxim S=\maxim T\cap S=\mathfrak nT\cap S=\mathfrak n$, showing that $f$ is also a scalar extension. Finally, assume
$f$ and $h$ are scalar extensions. Let
\begin{equation}
\label{eq:ffR}
\dots R^{b_2}\to R^{b_1}\to R\to R/\maxim\to 0
\end{equation}
be a free resolution of $R/\maxim$. Since $S$ is flat over $R$, tensoring yields a free resolution
\begin{equation}
\label{eq:ffS}
\dots S^{b_2}\to S^{b_1}\to S\to S/\maxim S\to 0.
\end{equation}
By assumption $S/\maxim S$ is the residue field $\lambda$ of $S$. Therefore,  $\tor S\bullet T \lambda $ can be calculated  
 as the   homology of
the complex
\begin{equation}
\label{eq:ffT}
\dots T^{b_2}\to T^{b_1}\to T\to T/\maxim T\to 0
\end{equation}
obtained from \eqref{eq:ffS} by the base change $S\to T$. However,
\eqref{eq:ffT} can also be obtained by tensoring \eqref{eq:ffR} over $R$ with $T$. Since $T$ is flat over $R$,
the sequence \eqref{eq:ffT} is exact, whence, in particular,  $\tor S1T \lambda =0$. By the local flatness
criterion,  $T$ is flat over
$S$. Since $\mathfrak n =\maxim S$ and $\pr=\maxim T$, we get   $\pr=\mathfrak nT$, showing that   $g$, too, is a scalar
extension.
\end{proof}

Three important examples of scalar extensions are given by the following
proposition.    

\begin{proposition}\label{P:ulsepsc} 
Let  $R$  be  a Noetherian local ring.
\begin{enumerate}
\item\label{i:compsc}  The natural map $R\to \complet R$ is a scalar extension.  
\item\label{i:et} Any etale map is a scalar extension.
\item\label{i:ulsepsc} The natural map $R\to \ulsep R$ is a scalar
extension, where $\ulsep R$ is any \c{}power of $R$.
\end{enumerate}
\end{proposition}
\begin{proof} 
The first two assertions  are well-known, so remains to show the
last. Let $\maxim$ be the maximal ideal of $R$. It is easy to show that $\maxim
\ulsep R$ is the maximal ideal of $\ulsep R$. So remains to prove  that $R\to
\ulsep R$ is flat.  Since $\ulsep R$ is complete, and in fact equal to the catapower of
$\complet R$, we may assume without loss of generality  that $R$ is already complete. In particular, $R$ is a homomorphic image of a regular local ring and if we prove the 
corresponding result for this regular local ring, then base change yields the
desired result by Lemma~\ref{L:ulsepbc}. Therefore, we may
moreover assume that $R$ is regular.  Since $\maxim \ulsep R$ is the maximal
ideal of $\ulsep R$ and since $\ulsep R$ is also regular by
\cite[Corollary 5.15]{SchFinEmb}, of the same dimension as $R$, the
flatness of $R\to \ulsep R$ then follows from \cite[Theorem 23.1]{Mats}. 
\end{proof}

In fact, \ref{P:ulsepsc}\eqref{i:et} has the following converse: if $R\to S$ is essentially of finite type
inducing a finite separable extension on the residue fields, then  $R\to S$ is a
scalar extension \iff\ it is etale. In this sense, scalar extensions are generalizations of etale
 maps (whence the alternative terminology `formally etale' for them).  This shows already that classification up to scalar extension is a
  reasonable and interesting problem. To gather further support for this
claim,  we will now explore how closely related scalar extensions are  to isomorphisms.
An important observation in that direction, one we will use several times
below, is that a
scalar extension of complete Noetherian local rings inducing an isomorphism on
their residue fields is itself an isomorphism; see \cite[Theorem 8.4]{Mats}. 
Hence it is of interest to generate scalar extensions $R\to S$ with $S$
complete. We will see that there exists a canonical choice over any field.

\subsection*{Completions along a residual extension} 
Let $(R,\maxim)$ be a Noetherian
local ring with residue field $\kappa$, and let $\lambda$ be a field extension of $\kappa$. The
\emph{completion of $R$ along $\lambda $} is the   (unique) local $R$-algebra $\scal
R \lambda $ solving the following universal problem: given an arbitrary  Noetherian
local $R$-algebra $S$ with residue field $\lambda $, if $S$ is complete, then there exists
a unique local $R$-algebra \homo\ $\scal R\lambda \to S$. When $\kappa=\lambda$,  we recover the
usual completion $\scal R\kappa=\complet R$ of $R$.  Here and
elsewhere, we say that there is a \emph{unique} \homo\ with certain properties, when we actually mean that 
there exists a unique \homo\ \emph{up to isomorphism}; this is consistent with
our practice of identifying two local rings when they are isomorphic. 

\subsubsection*{Proof of the existence of a completion along $\lambda $}
 We have to treat the equal
and mixed \ch\ cases separately. Assume first that $R$ has equal \ch\ (this case
is also discussed in \cite[(6.3)]{HHFreg}).
By Cohen's structure theorems, there exists an embedding $\kappa\to
\complet R$. Let $\scal R\lambda$ be the $\maxim(\complet R\tensor_\kappa\lambda)$-adic completion of $\complet
R\tensor_\kappa\lambda$. To see that this is a completion along $\lambda$, let $S$ be a
Noetherian local $R$-algebra with residue field $\lambda $ and assume $S$ is complete.
By the universal property of ordinary
completions, we get a unique \homo\ $\complet R\to S$. Since $S$ is complete,
we can find an embedding $\lambda\to S$ which agrees on the subfield $\kappa$ of $\lambda $ with
the composition $\kappa\to \complet R\to S$. By the
 universal  property of tensor products, the two maps $\complet R\to S$ and
$\lambda\to S$ combine  to a unique local $R$-algebra \homo\
$\complet R\tensor_\kappa \lambda\to S$, and using once more the universal property of
completion, this then yields a unique $R$-algebra \homo\ $\scal R \lambda\to S$. 

In the mixed \ch\ case, coefficient fields no longer exist and we now proceed as
follows. Let $V$ be the (unique)  complete $p$-ring  over $\kappa$,
where $p$ is the \ch\ of $\kappa$. We first  define the completion
of $V$ along $\lambda $, that is to say, $\scal V \lambda $, as the unique
complete $p$-ring over $\lambda $. That the latter satisfies the universal property of a completion along
$\lambda $ is proven in \cite[Theorem 29.2]{Mats}.  To define $\scal R \lambda $, let
$S$ be any Noetherian
local $R$-algebra with residue field $\lambda $ extending $\kappa$,  and assume $S$ is
complete. As before, we have a unique local $R$-algebra \homo\ $\complet R\to S$. By
 Cohen's structure theorems, there exists a commutative diagram of  local
\homo{s} 
\commdiagram[sc] V{}{\scal V \lambda}{} {} {\complet R}{} {S.}
By the universal property of
tensor products, we get a unique $R$-algebra \homo\ $\complet R\tensor_V\scal V \lambda\to
S$. Define $\scal R \lambda $ now as  the $\maxim(\complet R\tensor_V\scal V \lambda)$-adic completion of $\complet
R\tensor_V\scal V \lambda $, so that we get a unique local $R$-algebra \homo\ $\scal R \lambda\to S$, as
required.
\qed

\begin{corollary}\label{C:scal}
For every Noetherian local ring $R$ and every extension field $\lambda $ of its
residue field, $\scal R \lambda $, the  completion of $R$ along $\lambda $, exists and is
unique. For every ideal $I$ in $R$, the completion of $R/I$ along $\lambda $ is equal
to  $\scal R \lambda/I\scal R \lambda $. 

Moreover, the natural map $R\to \scal R \lambda $ is a scalar extension over $\lambda $.
\end{corollary}
\begin{proof}
Existence was proven   above; uniqueness then follows formally from being a
solution to a universal problem.
To prove the second assertion, assume $R/I\to S$ is a
local \homo\ with $S$ a complete Noetherian local ring with residue field $\lambda $.
The composition $R\to R/I\to S$ yields by definition a unique local $R$-algebra
\homo\ $\scal R \lambda\to S$. Since $IS=0$, the latter \homo\ factors through $\scal R \lambda/I\scal
R \lambda $, showing that $\scal R \lambda/I\scal R \lambda $ satisfies the universal property of
completions along $\lambda $.
As for the last assertion, in the equal
\ch\
case, the base change $\complet R\to \complet R\tensor_\kappa \lambda $ of $\kappa\sub \lambda $  is
faithfully flat. Since completion is exact,  each map in  
\[
R\to \complet R\to \complet R\tensor_\kappa \lambda\to \scal R \lambda
\]
is faithfully flat, whence so is their composition.  In the mixed \ch\ case,
$\scal V \lambda $ is
torsion-free whence flat over $V$. Hence by the same argument as in the equal
\ch\ case, the composite map 
\[
R\to \complet R\to \complet R\tensor_V\scal V \lambda\to \scal R \lambda
\]
 is faithfully flat. By our second assertion,  $\scal
R \lambda/\maxim\scal R \lambda $ is the completion of $R/\maxim\iso \kappa$ along $\lambda $ in either
\ch.
In other words,  $\scal R \lambda/\maxim\scal R \lambda\iso \lambda $ and hence
in particular, $\maxim\scal R \lambda $ is the maximal ideal of $\scal R \lambda $. This
proves that $R\to \scal R \lambda $ is a scalar extension.
\end{proof}

\begin{proposition}\label{P:compscal} 
Let $R\to S$ be a scalar extension over $\lambda $.
If $S$ is complete, then $S\iso \scal R \lambda $.
\end{proposition}
\begin{proof} 
By the universal property, we have a local $R$-algebra \homo\ $\scal R \lambda\to S$. It follows from
\cite[Theorem 8.4]{Mats} that $\scal R \lambda\to S$ is   surjective. Since $R\to
\scal R \lambda $ and $R\to S$ are scalar extensions by Corollary~\ref{C:scal}  and 
  by assumption respectively, $\scal R \lambda\to S$ is faithfully flat 
by Proposition~\ref{P:scalfact}, whence injective.
\end{proof}

\begin{corollary}[(Lifting of scalar extensions)]\label{C:liftscal}
Let $R\to S$ be a scalar extension with  $S$ complete. If $R$ is the
homomorphic image of a Noetherian local ring $A$, then there exists a scalar
extension $A\to B$ whose base change is $R\to S$, that is to say,
$S=B\tensor_AR$.
\end{corollary}
\begin{proof}
We leave it to the reader to verify that, after taking completions, we may assume
that also $A$ and $R$ are complete. By Cohen's structure theorems, $A$ and $R$ are
the homomorphic images of $\pow VX$ modulo some ideals $J\sub I$ respectively,
where $V$ is either their common residue
field or otherwise a complete $p$-ring over that residue field, and where
$X$ is a finite tuple of indeterminates. Moreover, $S\iso \scal R \lambda $ by
Proposition~\ref{P:compscal}, where
$\lambda $ is the residue field of $S$. In particular, $S\iso
\pow {\scal V \lambda}X/I\pow{\scal V \lambda}X$. Hence putting $B:=\pow{\scal
V \lambda}X/J\pow{\scal V \lambda}X$ yields a scalar extension $A\to B$ with $A/IA=R\to
B/IB=S$, as required.
\end{proof}

The following result is a sharpening of \cite[Theorem 2.4]{OlbSaySha}.

\begin{corollary}\label{C:ulsepscal}
Let $R$ be a Noetherian local ring with residue field $\kappa$. If $\ul \kappa$ is
the ultrapower of $\kappa$, then the \c{}power $\ulsep R$ of $R$ is equal to the 
completion $\scal R{\ul \kappa}$ along $\ul \kappa$.
\end{corollary}
\begin{proof}
By Lemma~\ref{L:ulsepbc}, the residue field of $\ulsep R$ is $\ul \kappa$. Since $R\to \ulsep R$ is a scalar extension
by \ref{P:ulsepsc}\ref{i:ulsepsc}, and since $\ulsep R$ is complete, $\ulsep R\iso \scal
R{\ul \kappa}$ by Proposition~\ref{P:compscal}.
\end{proof}

\begin{corollary}\label{C:finscal}
Let $R\to S$ be a finite local \homo\ inducing a trivial extension on the residue
fields. For every extension $\lambda $ of this common residue field, $\scal
S \lambda\iso\scal R \lambda\tensor_RS$.
\end{corollary}
\begin{proof}
The base change $S\to \scal R \lambda\tensor_RS$ is faithfully flat. Let $\maxim$ and
$\mathfrak n$ be the maximal ideals of $R$ and $S$ respectively. Since
\[
(\scal R \lambda\tensor_RS)/\mathfrak n(\scal R \lambda\tensor_RS)\iso(\scal R \lambda/\maxim\scal 
R \lambda)\tensor_{R/\maxim} (S/\mathfrak n)\iso \lambda\tensor_\kappa\kappa= \lambda 
\]
the ideal $\mathfrak n(\scal R \lambda\tensor_RS)$ is a maximal ideal. Since the base
change $\scal R \lambda\to \scal R\lambda\tensor_RS$ is finite with trivial residue field
extension and since $\scal R \lambda $ is complete whence Henselian,
$\scal R \lambda\tensor_RS$ is a complete local ring. Hence we showed that $S\to  \scal
R \lambda\tensor_RS$ is
a scalar extension and since the latter ring is complete with residue field
equal to $\lambda $, it is isomorphic to $\scal S \lambda $ by Proposition~\ref{P:compscal}.
\end{proof}

\begin{corollary}\label{C:scalreg}
Suppose $R$ is  an    excellent local ring. If
$R\to S$ is a scalar extension inducing a separable extension on the residue
fields, then $R\to S$ is a regular
\homo.
\end{corollary}
\begin{proof}
By \cite[Theorem 28.10]{Mats}, the scalar extension $R\to S$ is formally smooth,
since it is unramified and the residue field extension is separable. The
assertion now follows from a result by Andr\'e in \cite{And} (see also \cite[p.
260]{Mats}). 
\end{proof}

In fact, with aid of  Proposition~\ref{P:compscal}, Corollary~\ref{C:liftscal}
and Cohen's structure theorems, one reduces to proving that $\pow VX\to \pow{\scal
V \lambda}X$ is regular, where $V$ is either a field or a complete $p$-ring, and where
$\lambda $ is a separable extension of the residue field of $V$. This approach circumvents the use of Andr\'e's deep result.

\begin{definition}\label{D:irr}
A Noetherian local ring $R$ is called \emph{analytically irreducible}, if
$\complet R$ is a domain; it is called \emph{absolutely analytically
irreducible}, if $\scal R{\ac \kappa}$ is a domain, where $\ac \kappa$ is the
algebraic closure of the residue field $\kappa$ of $R$; and it is called  
\emph{universally irreducible},
if any scalar extension of $R$ is a domain.
\end{definition}

\begin{corollary}
If $R$ is an excellent normal local domain with perfect residue field, then $R$
is universally irreducible.
\end{corollary}
\begin{proof}
Let  $S$ be a scalar extension of $R$. By Corollary~\ref{C:scalreg}, the map
$R\to S$ is regular and hence $S$ is again normal by \cite[Theorem 32.2]{Mats},
whence a domain. 
\end{proof}

\begin{proposition} 
A Noetherian local ring is absolutely analytically irreducible \iff\ it is universally
irreducible.
\end{proposition}
\begin{proof}
Since we will make no essential use of this result, we only give a sketch of a
proof. One direction is obvious. For the other, we may reduce to the case that $R$ is a complete Noetherian
local domain with algebraically closed residue field $\kappa$. We need to show that
$\scal R \lambda $ is a domain, where $\lambda $ is an arbitrary extension field of $\kappa$. By
Cohen's structure theorems, there exists a  finite extension $S:=\pow VX\sub R$,
where $V$ is either $\kappa$ or the complete $p$-ring over $\kappa$, and $X$ is a tuple of
indeterminates. Write $R=\pol SY/\pr$ for some finite tuple of indeterminates
$Y$, so that $\pr$ is in particular a prime ideal.
Since the fraction field of $\scal S \lambda =\pow{\scal V \lambda}X$ is a regular extension
of the fraction field of $S=\pow VX$,  the same argument as in the proof
of \cite[Lemma 5.21]{SchAsc} then shows that $\pr\pol{\scal S \lambda}Y$ is a prime
ideal. Hence we are done, since $\scal R \lambda =\pol{\scal S \lambda}Y/\pr \pol{\scal S \lambda}Y$ by
Corollary~\ref{C:finscal}.
\end{proof}

We are ready to formulate a flatness criterion generalizing
\cite[Theorem 8]{Kol}; we prove a slightly stronger version than the one quoted in the introduction.

\begin{theorem}\label{T:ff}
Let $R\to S$  be a local \homo\ of Noetherian
local rings. Assume $R$ is universally irreducible, e.g., an excellent
normal local domain with perfect residue field, or a complete local
domain with algebraically closed residue field. If $R\to S$ is unramified and
$\op{dim}(R)=\op{dim}(S)$, then $R\to S$ is faithfully flat, whence a scalar
extension.
\end{theorem}
\begin{proof} 
Recall that $(R,\maxim)\to(S,\mathfrak n)$ being unramified means that $\mathfrak
n=\maxim S$. It suffices to prove the assertion under the additional
assumption that both
$R$ and $S$ are complete. Indeed, if $R\to S$ is arbitrary, then $\complet
R\to \complet S$ satisfies again the hypotheses of the theorem and therefore
would be faithfully flat. By an easy descent argument, $R\to S$ is
then also faithfully flat.

So assume $R$ and $S$ are complete and let $\lambda $ be the residue field of $S$.  By
assumption,  $\scal R \lambda $ is   a domain, of the same dimension as $R$.  By the
universal property of the completion along $\lambda $, we get a local $R$-algebra
\homo\ $\scal R \lambda\to S$. By
\cite[Theorem 8.4]{Mats}, this \homo\ is surjective. It is also injective, since 
$\scal R \lambda $ and $S$ have the same dimension and $\scal R \lambda $ is a domain.  Hence
$\scal R \lambda\iso S$, so that $R\to S$ is a scalar extension.
\end{proof}

\section{Similarity relation}\label{s:sim}

 Next, we introduce an equivalence relation on the class of Noetherian local
rings which, although
coarser than the isomorphism relation, preserves most local singularity
properties (see for instance Theorem~\ref{T:transsim} below). Namely, we say that
two Noetherian local rings $R$ and $S$  are
\emph{similar}, denoted $R\approx S$, if they admit a common scalar extension.  Let $T$ be this common scalar extension.
Its completion is again a scalar extension and by Proposition~\ref{P:compscal}, it is
therefore isomorphic to both $\scal R \lambda $ and $\scal S \lambda $, where $\lambda $ is the residue
field of $T$. In other words, we showed that $R\approx S$
\iff\  $\scal R \lambda\iso\scal S \lambda $ for some  sufficiently large common extension $\lambda $
of their respective residue fields. It follows easily from this that $\approx$ is an equivalence
relation. The collection of all local rings similar to a given Noetherian local ring $R$ is
called the \emph{similarity class} of  $R$ and is denoted $[R]$. Immediately from the results in \cite[\S23]{Mats} and
\cite[Proposition 9.3]{SchCI} (where the notion of a \emph{singularity defect} is
introduced), we get:

\begin{theorem}\label{T:transsim}
If two Noetherian local rings   are similar, then they have the same dimension,
depth and Hilbert series, and one is regular (respectively, \CM, Gorenstein,
complete intersection) \iff\ the other is. More generally, any two similar local rings
have the same singularity defects.\qed
\end{theorem}

Using Corollary~\ref{C:scalreg}, other properties, such as being reduced or normal, are also invariant under the
similarity relation, provided the rings are excellent with perfect residue
field. Note that being a domain is not preserved under the similarity
relation, necessitating definitions~\ref{D:irr}.

\begin{proposition}\label{P:ulsepsim} 
Any two \c{}powers of a Noetherian local ring, or more generally, any
two  Noetherian local rings which are elementary equivalent,   are similar.

More generally, let   $\seq Rw$ and $\seq Sw$ be
sequences of Noetherian local rings of    embedding dimension at most $d$. If
almost each $\seq Rw$ is similar to $\seq Sw$,   then the respective \c{}products
$\ulsep R$ and $\ulsep S$ are also similar.
\end{proposition}
\begin{proof}
Suppose $R$ and $S$ are elementary equivalent Noetherian local rings. By the
Keisler-Shelah theorem (see \cite[Theorem 9.5.7]{Hod}), some ultrapower of $R$
and $S$ are isomorphic, whence so are their corresponding \c{}powers
(strictly speaking, the underlying index set will in general no longer  be
countable, so that we have to make some   minor modifications alluded to in
footnote~\eqref{f};   details are left to the reader).
By Proposition~\ref{P:ulsepsc}, these are scalar extensions of $R$ and $S$
respectively, proving the first assertion.

To prove the second assertion, we may without loss of generality assume that all
rings are complete. By our discussion above, we may further reduce to
the case that $\seq Sw$ is a scalar extension of $\seq Rw$. Since $\seq Rw$ is a
homomorphic
image of a $d$-dimensional regular local ring by Cohen's structure theorems, and
since the property we seek to prove is preserved under homomorphic images by
Lemma~\ref{L:ulsepbc} and Corollary~\ref{C:scal}, we may moreover assume by
Corollary~\ref{C:liftscal} that
each $\seq Rw$ is regular, of dimension $d$. By Theorem~\ref{T:transsim},
almost all $\seq Sw$ are then also regular of dimension $d$. By \cite[Corollary 5.15]{SchFinEmb}, the \c{}products $\ulsep R$ and
$\ulsep S$ are therefore again $d$-dimensional regular local rings.    The induced  \homo\ $\ulsep R\to
\ulsep S$   is unramified by Lemma~\ref{L:ulsepbc}. Hence, it is
faithfully flat by \cite[Theorem 23.1]{Mats}, whence a scalar extension, as we wanted to show. 
\end{proof}

We denote the collection of all 
similarity classes of  Noetherian local rings   by
$\simloc$.  
Although the class of Noetherian local rings is not a set, we do no longer have
this complication for its quotient:

\begin{proposition}\label{P:set}
 The quotient $\simloc$ is a set.
\end{proposition}
\begin{proof} 
Let $[R]$ be a similarity class   and let $\kappa$ be the
residue field of $R$. Since $R\approx
\complet R$, we may assume that $R$ is complete, whence,  by   Cohen's structure
theorems, the homomorphic image of  $S:=\pow VX$ with $V$ either
equal to $\kappa$ or to the   complete $p$-ring over  $\kappa$, and with 
$X$ a finite tuple of indeterminates. Suppose $R=S/I$ with $I=\rij fsS$. We may
choose a  
subring $W$ of $V$ of size at most the continuum so that it contains all coefficients of the $f_i$  
and so that $W$ is again  a
field or a  complete  $p$-ring. Let $T:=\pow{W}X$ and $J:=\rij fsT$, 
so that $S\iso\scal T\kappa$ and $I=JS$.
Hence, by base change, $R\iso S/I$ is a scalar extension of $T/J$, showing that
$T/J\approx R$. In conclusion, we showed that every similarity class
contains
a ring of size at most the continuum, and therefore $\simloc$ is a set.
\end{proof}

\section{\Deform\ metric}

Our next goal is to define a metric on the space $\simloc$. We will first
define a semi-metric on the space of all Noetherian local rings.

\subsection*{\Deform\ semi-metric}
Let $(R,\maxim)$ be a Noetherian local ring. The \emph{$n$-th \infdef}
of $R$ (also called the \emph{$n$-th infinitesimal neighborhood}) is by definition the (Artinian) residue ring $R/\maxim^n$ and will be
denoted $\app Rn$. Recall that the ($\maxim$-adic) completion $\complet R$
of $R$ is the inverse limit of all $n$-th \infdef{s} of $R$, and that $\app
Rn\iso \app{\complet R}n$.  We define a semi-metric on the class of all  Noetherian local rings, called
the
\emph{\deform\ metric}, as follows. Given 
two  Noetherian local rings $R$ and $S$, let $d(R,S)$ be  the  infimum of the
numbers $2^{-n}$ for which  $\app Rn\iso \app Sn$. In words, the distance between
two  local rings  is at most $2^{-n}$ if their $n$-th \infdef{s}
agree.  One easily verifies that this distance function satisfies all the axioms
of a metric, except that two
distinct elements can be at distance zero, so that $ d(\cdot,\cdot)$ is only a  
semi-metric. It is an interesting problem to determine all local rings that
are $d$-equivalent to a
given local ring; a partial answer is provided in \cite{vdDComp}. It is clear that any two Noetherian local
rings with the same completion have this property. For our purposes, the
following partial solution to this question suffices:

\begin{proposition}\label{P:distzero}
Given two Noetherian local rings $R$ and $S$, if $R\sim_dS$, that is to say, if $
d(R,S)=0$, then $R\approx S$.
\end{proposition}
\begin{proof}
By definition, there exists  for each $n$ an Artinian local ring $\seq Tn$
isomorphic to both  $\app Rn$ and $\app Sn$. Let  $\ul R$ and $\ulsep R$ be the respective ultrapower and \c{}power
of $R$, and let   $\ul T$ and $\ulsep T$ be  the
respective ultraproduct and \c{}product of  the $\seq Tn$. Taking ultraproducts of the
surjections $R\to \seq Tn$ yields a surjection $\ul R\to\ul T$
whence a surjection $\ulsep R\to\ulsep T$. Let $r\in\ul R$ be an
element whose image in $\ulsep R$ lies in the kernel of $\ulsep R\to\ulsep
T$, that is to say, $r\in\inf{\ul T}$. Let $\seq
rn$ be elements in $R$ with ultraproduct equal to $r$. Fix some $N$. Since 
$r\in\maxim^N\ul T$, \los\ yields  $\seq rn\in\maxim^NT_n$ for almost
all $n$. For $n\geq
N$, this implies $\seq rn\in\maxim^N$ and hence by \los, $r\in\maxim^N\ul R$.
Since this holds
for all $N$, the image of $r$ in $\ulsep R$ is zero, showing that $\ulsep
R\to\ulsep T$ is an isomorphism. Applying the same argument to the
\c{}power $\ulsep S$ of $S$, we also
get $\ulsep S\iso\ulsep T$ and hence $\ulsep R\iso\ulsep S$.
Therefore, $R\approx S$ by Proposition~\ref{P:ulsepsc}\eqref{i:ulsepsc}.
\end{proof}

The \deform\ semi-metric is  non-archimedean, and hence the induced topology, 
called the \emph{\deform\
topology}, is totally disconnected. By convention, the zero-th \infdef\ of a ring is zero (since we think of
$\maxim^0$ as the unit ideal). It follows
that the distance between any two local rings is at most one, that is to say,
$d$ is bounded. Immediately from the definitions we also get:

\begin{lemma}\label{L:embdim} 
If $ d(R,S)<1$, then $R$ and $S$ have the same residue field; if $ 
d(R,S)<1/2$, then $R$ and $S$ have the same embedding dimension.\qed
\end{lemma}

In particular, if, in this metric, $\seq Rw$ is a Cauchy sequence of Noetherian local rings,
then    almost of all $\seq Rw$ have the same residue
field, called the \emph{residue field}  of the sequence, and the same embedding
dimension. By the above discussion, the \c{}product $\ulsep R$ is therefore a complete Noetherian
local ring. By Lemma~\ref{L:embdim}, the embedding dimension is a continuous map
onto the discrete space $\zet$.  This is no longer
true for dimension: for instance $R:=\pow kX$ and $R_n:=R/X^nR$ lie at distance $2^{-n}$, yet their dimensions are not
the same.  One can show, however, that dimension is upper-semicontinuous.

By an  (open) \emph{ball} $\ball$ with
\emph{center} $R$ and \emph{radius} $0< \delta\leq 1$, we mean the collection of
all Noetherian local rings $S$  such that $  d(R,S)<\delta$.  Since
the metric is non-archimedean, any member of a 
 ball is its center and every    ball is both open and closed in the \deform\
topology, that is to say, is a \emph{clopen}.
 Because the distance function only takes discrete values (the powers of $1/2$), any two radii
 which lie between two consecutive powers of $1/2$ yield the same ball. Therefore,
by the \emph{radius} of a ball $\ball$, we mean twice
the largest distance between two members of $\ball$; this is always a 
power of $1/2$. (We need to take twice the distance since we used a strict
inequality in the definition of a ball.)

A \emph{unit ball} is a ball $\ball$ with radius $1$ and hence consists of all local rings with the same residue field.
We call this common residue field the \emph{residue field} of $\ball$.  This gives a one-one correspondence between unit
balls and fields.   More generally, to every ball $\ball$, we associate an
Artinian local ring $R_{\ball}$, called the
\emph{residue ring} of $\ball$, given as the unique local ring such that
$\app Rn\iso R_{\ball}$, for all $R\in\ball$, where $2^{-n+1}$ is the radius of $\ball$. Note that
$R_\ball$ is a center of $\ball$ and, moreover,
the radius of $\ball$ is determined by $R_{\ball}$: it is equal
 to $2^{-n+1}$ where $n$ is the
 nilpotency index of $R_{\ball}$. In conclusion, there is a
one-one correspondence between balls 
$\ball$   and    Artinian local rings.

\begin{proposition}\label{P:card} 
Every ball   is a set.
\end{proposition}
\begin{proof} 
It suffices to prove this for a unit ball $\ball$. The result will
follow if we can show that there is a
cardinal number so that every member of $\ball$ has size at most this cardinal. Let $\kappa$ be the residue
field of $\ball$ and let $R\in\ball$. Since the cardinality  of a Noetherian local ring is at most the cardinality of
its completion, we may assume that $R$ is complete.  By   Cohen's structure
theorems, $R$ is a homomorphic image of $\pow VX$, where $X$ is a finite tuple of indeterminates and  
$V$ is equal to $\kappa $ in the equal \ch\ case, and equal to the complete
$p$-ring over $\kappa $ in the mixed \ch\ case.    It is clear that in either case,
the cardinality of $\pow VX$ is bounded in
terms of the cardinality of $\kappa$, whence so is its homomorphic image $R$.
\end{proof}

Note that each ball $\ball$ is infinite: if $R_\ball$ is its residue ring, then
the latter is of the form $S/I$, where $(S,\mathfrak n)$ is a power series ring $\pow
VX$. If $n$ is the nilpotency index
of $R_\ball$, then $S/J\in\ball$ for any ideal $J\sub S$ such that
$J+\mathfrak n^n=I$.

\begin{corollary}\label{C:isomball}
Let $\kappa\subseteq \lambda $ be an extension of fields and let $\ball_\kappa$
and $\ball_\lambda$ be the unique unit balls with residue field $\kappa$ and $\lambda $, respectively. The map sending a ring in $\ball_\kappa$ to its
completion along $\lambda $ is an isometry $\ball_\kappa\to\ball_\lambda$.
\end{corollary}
\begin{proof} 
Take $R,S\in\ball_\kappa$. Clearly, the completions $\scal R \lambda $ and $\scal S \lambda $ along
$\lambda $ belong both to $\ball_ \lambda $.  Suppose $  d(R,S)\leq 2^{-n}$,
that is to say, their $n$-th \infdef{s} $\app Rn$
and $\app Sn$ are isomorphic. By Corollary~\ref{C:scal}, the completions of $\app
Rn$ and $\app Sn$ along $\lambda $ are respectively $\app{\scal R \lambda}n$ and $\app{\scal
S \lambda}n$, and
therefore are  isomorphic, showing that $  d(\scal R \lambda,\scal S \lambda)\leq 2^{-n}$.
\end{proof}

\begin{proposition}\label{P:lim} 
If $\mathbf r$ and  $\mathbf s$ are Cauchy sequences  of Noetherian local
rings, say, $\tuple r(w):=\seq Rw$ and $\tuple s(w):=\seq Sw$, with respective \c{}products $\ulsep R$ and $\ulsep S$,
then  $ 
d(\ulsep R,\ulsep S)\leq  d(\mathbf r,\mathbf s)$.    In
particular,  if $\mathbf r\sim_d\mathbf s$, then  $\ulsep R \approx
\ulsep S $. 
\end{proposition}
\begin{proof} 
The last assertion is immediate by the first and Proposition~\ref{P:distzero}.
Suppose $  d(\mathbf r,\mathbf s)\leq 2^{-n}$. This means that
for some $w_0$ and all $w>w_0$, we have $\app{\seq Rw}n\iso \app{\seq Sw}n$. By Lemma~\ref{L:ulsepbc}, the $n$-th \infdef{s} $\app{\ulsep R}n$ and $\app{\ulsep S}n$ are isomorphic, showing that  
$  d(\ulsep R,\ulsep S)\leq 2^{-n}$. 
\end{proof}

The next result shows that \c{}products act as  limits up to similarity.
To formulate it, we extend our previous notation: let $\mathbf
r$ be a sequence of Noetherian local rings with the same residue
field $\kappa$ (e.g., a Cauchy sequence) and let $\lambda $ be an extension field of $\kappa$.
Then we let $\scal {\mathbf r} \lambda $ denote the sequence of rings obtained
by taking the completions along $\lambda $ of all members of $\mathbf r$, that is to
say, $\scal{\mathbf r} \lambda(w):=\scal{\seq Rw} \lambda $, if $\tuple r(w)=\seq Rw$.

\begin{theorem}\label{T:lim} 
Let $\mathbf r$ be a Cauchy sequence of  Noetherian local rings    with   
residue field $\kappa$. Let $\lambda $ be any extension field of the ultrapower $\ul \kappa$ of $\kappa$.
Then $\scal {\mathbf r} \lambda $
is a Cauchy sequence converging to $\scal{\ulsep R} \lambda $. In
particular, $\ulsep R$ is a limit of $\scal {\mathbf r}{\ul \kappa}$.
\end{theorem}
\begin{proof}
Let $\seq Rw:=\tuple r(w)$. Fix $n$ and choose $w(n)$ so that all $\app{\seq Rw}n$ for $w\geq w(n)$ are
isomorphic, say, to $T$.   By Lemma~\ref{L:ulsepbc}, the $n$-th \infdef\
$\app{\ulsep R}n$
is isomorphic to the \c{}power $\ulsep {T}$; the latter is
isomorphic to $\scal T{\ul \kappa}$ by Corollary~\ref{C:ulsepscal}; and this in
turn is isomorphic to $\app{\big(\scal {\seq Rw}{\ul \kappa} \big)}n$, for all
$w\geq w(n)$ by Corollary~\ref{C:scal}. In
summary, we showed that 
$$
  d(\scal {\seq Rw}{\ul \kappa},\ulsep R)\leq
2^{-n},
$$
 for all $w\geq w(n)$. By Corollary~\ref{C:isomball}, taking
completions along $\lambda $ yields
$  d(\scal {\seq Rw} \lambda,\scal{\ulsep R} \lambda)\leq 2^{-n}$, for all
$w\geq w(n)$.  Since this holds for all $n$, the assertion follows.
\end{proof}

\section{Similarity space}

We are ready to define a metric on the similarity space $\simloc$.   For two similarity classes
$[R]$ and $[S]$,  let $  d([R],[S])$ be equal to the infimum of all $ 
d(R',S')$ with $R'\approx R$ and
$S'\approx S$. Alternatively, recall that for a   semi-metric space $(\Sigma,d)$,
the distance $  d(U,V)$ between two subclasses $U$ and
$V$ is defined to be the infimum of all $  d(x,y)$ with $x\in U$ and $y\in V$; hence $ 
d([R],[S])$ is just the distance between $[R]$ and $[S]$ viewed as  subclasses.

\begin{lemma}\label{L:app} 
For any two Noetherian local rings $R$ and
$S$ and  for any $n\in\nat$, we have $  d([R],[S])\leq 2^{-n}$
\iff\ $\app Rn\approx \app Sn$.
\end{lemma}
\begin{proof} 
Suppose $  d([R],[S])\leq 2^{-n}$ and choose $R'\approx R$  and $S'\approx S$ so that
$  d(R',S')\leq 2^{-n}$. In other words,  $\app {R'}n\iso \app {S'}n$ and
therefore,  $\app Rn\approx \app Sn$  by  Corollary~\ref{C:scal}.
Conversely, assume $\app Rn\approx
\app Sn$ and let $T$ be a common scalar extension of $\app Rn$ and $\app
Sn$. Let $\lambda $ be the residue field of $T$. By
Corollary~\ref{C:scal}, the  $n$-th \infdef{s}
of $\scal R \lambda $ and $\scal S \lambda $ are equal to $T$.  In other words, $ 
d(\scal R \lambda,\scal S \lambda)\leq 2^{-n}$. Since
$  d([R],[S])$ is defined as an infimum,  it is at most  $2^{-n}$.
\end{proof}

 \begin{corollary}\label{C:metric}
 The quotient $(\simloc, d)$  is a metric space.
\end{corollary}
\begin{proof} 
Suppose $  d([R],[S])=0$. By Lemma~\ref{L:app}, the $n$-th
\infdef{s} $\app Rn$ and $\app Sn$ of $R$ and $S$ are similar for all $n$. Hence
there exists a common scalar extension $T_n$ of $\app Rn$ and $\app Sn$. 
We may inductively choose $T_{n+1}$ to have a residue field containing the residue field of $T_n$ by
Corollary~\ref{C:isomball}, since scalar extensions can only make the distance
smaller. Let $\lambda $ be the union of all these residue fields. By
Corollary~\ref{C:scal}, the $n$-th \infdef{s} of $\scal R \lambda $ and $\scal S \lambda $
are equal to  $\scal {T_n} \lambda $. Since
this holds for all $n$, we showed that $  d(\scal R \lambda,\scal S \lambda)=0$. By
Proposition~\ref{P:distzero}, we get  $\scal R \lambda\approx \scal S \lambda $ and hence
$[R]=[\scal R \lambda]=[ \scal S \lambda]=[S]$.
\end{proof}

It follows from Theorem~\ref{T:lim} that given a Cauchy sequence $\mathbf r$
of  Noetherian local rings,
  the sequence  $\scal{\mathbf r}{\ul \kappa}$   has a limit, where $\ul
\kappa$ is the  ultrapower of the  residue field of $\mathbf r$. Since the
corresponding members
of $\mathbf r$ and $\scal{\mathbf r}{\ul \kappa}$ are similar, we showed that every
Cauchy sequence becomes convergent after
replacing each of its components by an appropriately chosen similar ring. Therefore, the next result should not come as
a surprise:

 \begin{theorem}\label{T:simcomp} 
 The metric space $\simloc$ is complete.
\end{theorem}
\begin{proof} 
We will define an isometry $\complet \imath\colon \complet\simloc\to
\simloc$ as follows. We start with defining a map $i\colon \cau\simloc \to
\simloc$. Let ${\mathbf r}$ be a Cauchy sequence in
$\simloc$. For each $w$, let $\seq Rw$ be a representative in the similarity class $\tuple r(w)$, and let 
 $\ulsep R$ be their \c{}product. Note that $\ulsep R$ is a complete Noetherian local ring since almost all $\seq Rw$ have the
same embedding dimension.
Define $i({\mathbf r}):=[\ulsep R]$. By
Proposition~\ref{P:ulsepsim}, the map $i$ is well-defined, that is to say,
does not depend on the choice of representatives $\seq Rw$. 
Suppose ${\mathbf s}$ is a second
Cauchy sequence which is equivalent to ${\mathbf r}$ and let $\ulsep S$  be the \c{}product of the representatives 
$\seq Sw$ of each $\tuple s(w)$. For a fixed $n$, we have $  d([\seq Rw],[\seq Sw])\leq 2^{-n}$ for all sufficiently large $w$. By
Lemma~\ref{L:app}, the $n$-th \infdef{s} of $\seq Rw$ and $\seq Sw$ are
therefore similar, for all sufficiently large $w$. By
Proposition~\ref{P:ulsepsim}, then so are the $n$-th \infdef{s} of  $\ulsep R$ and $\ulsep S$,  so that  $  d([\ulsep R],[\ulsep S])\leq 2^{-n}$ by another application of Lemma~\ref{L:app}.
Since this holds for all $n$, Corollary~\ref{C:metric} yields $[\ulsep R]=[\ulsep S]$.
 By definition of completion, $i$ therefore factors through a map 
 \[
 \complet \imath\colon \complet\simloc\to \simloc.
 \]
  We leave it to the reader to check
that $\complet \imath$ preserves the metric. Note that $\complet \imath$
restricted to $\simloc$ is the identity, since a \c{}power is a scalar
extension by Proposition~\ref{P:ulsepsim}. Hence $\complet \imath$
must be surjective. To prove injectivity, assume $\mathbf r$ and $\mathbf s$
are Cauchy sequences of  Noetherian local rings  whose respective \c{}products
$\ulsep R$ and $\ulsep S$ are similar. Let $\lambda $ be a large
enough field extension so that 
\[
\scal{\ulsep R} \lambda\iso \scal{\ulsep S} \lambda.
\]
 By Theorem~\ref{T:lim}, the (component-wise) completion $\scal{\tuple
r} \lambda $ along $\lambda $ converges to $\scal{\ulsep R} \lambda $, and likewise $\scal{\mathbf s} \lambda $
 converges to $\scal{\ulsep S} \lambda $. Therefore, $\scal{\mathbf
r} \lambda $ and $\scal{\mathbf s} \lambda $, as they converge to the same limit, are equivalent, which 
proves that $\complet \imath$ is injective.
\end{proof}

We have the following  generalization of Proposition~\ref{P:ulsepsim}.

\begin{corollary}
 If $\seq Rw$ is a Cauchy sequence, then any two cataproducts  of $\seq Rw$
(with respect to different ultrafilters) are similar. In particular, if the common residue field $\kappa$ of the $\seq Rw$ is an   \acf, then the cataproduct of the $\seq Rw$ is, up to isomorphism, independent from the choice of ultrafilter.
\end{corollary}
\begin{proof}
According to the proof of Theorem~\ref{T:simcomp}, the similarity class of any
cataproduct $\ulsep R$ of the $\seq Rw$ is a limit of the sequence of similarity
classes $[\seq Rw]$, and therefore, is unique by Corollary~\ref{C:metric}.

Since any two   ultrapowers of $\kappa$ are algebraically closed and have the same (uncountable) cardinality, they are isomorphic by Leibnitz's theorem.  Since any two cataproducts of the $\seq Rw$ are similar by the first assertion, and are complete with isomorphic residue fields, they must be isomorphic by Proposition~\ref{P:compscal}.
\end{proof}

We introduce the following notation. Let $\mathbb S\sub \simloc$ be a subset, and let $d\geq 0$ and $e\geq 1$. We let $\mathbb S_d$ (respectively, $\mathbb S_{d,e}$) be the set of similarity classes  of Noetherian local rings in $\mathbb S$ having dimension $d$ (and  parameter degree $e$). Recall that the \emph{parameter degree} of $R$ is the minimal length of a residue ring $R/I$, where $I$ runs over all parameter ideals of  $R$. It is not hard to show that two similar rings with infinite residue field have the same parameter degree, and so we may speak of the parameter degree of a similarity class as the parameter degree of any of its members having infinite residue field.

\begin{corollary}
For each $d\geq 0$ and $e\geq 1$, the subset $\simloc_{d,e}\sub\simloc$ is closed.
\end{corollary}
\begin{proof}
It suffices to show that $\simloc_{d,e}$ is closed under limits. Hence let $\tuple r$ be a Cauchy sequence in $\simloc_{d,e}$, and choose representatives $\seq Rw$ in each $\tuple r(w)$, of dimension $d$ and parameter degree $e$. Let $\ulsep R$ be the cataproduct of the $\seq Rw$, so that its similarity class is  the limit of $\tuple r$ by Theorem~\ref{T:simcomp}. Since the cataproduct $\ulsep R$ has dimension $d$ and parameter degree $e$ by \cite[Theorem 5.22]{SchFinEmb},  the claim follows.
\end{proof}

We can now state and prove the main theorem of this paper:

 \begin{theorem}\label{T:smooth}
 The metric space $\simloc$ is a Polish space. In particular,
the similarity relation is smooth.
\end{theorem}
\begin{proof}
In view of Theorem~\ref{T:simcomp}, it remains to show that
$\simloc$ contains a countable dense subset. We already observed
that there is a one-one correspondence between balls  and
Artinian local rings, so that $\simloc_0$,  the similarity classes of Artinian local rings, form
a dense subset of  $\simloc$. Let $R$ be an Artinian local ring
with
residue field $\kappa$. By Cohen's structure theorems, $R$ is of the form $\pow VX/I$, where $V$ is either $\kappa $ or  
 the complete $p$-ring over $\kappa $, and where $X$ is  a tuple of indeterminates. 
Since $R$ is Artinian, it is in fact
finitely generated over $V$.  Hence, by an argument similar to the one  in the proof of
Proposition~\ref{P:set}, there
exists a finitely generated subfield $\kappa_0\sub \kappa$ and
an Artinian local ring $R_0$ with residue field $\kappa_0$, such that $R_0\approx R$.
Since there are only countably many finitely generated fields,
the collection of all these $R_0$ is again countable.
\end{proof}

\section{Variants}

A  first variant is simply obtained by working in the category   of all Noetherian local $Z$-algebras, for $Z$
some Noetherian ring, so that the morphisms are now given by local $Z$-algebra
\homo{s}. This leads to the notion of two $Z$-algebras being
\emph{$Z$-similar}, and the same argument shows that classifying Noetherian
local $Z$-algebras up to $Z$-similarity is again a smooth problem. 

We may also extend the definition to include modules. Namely, given an $R$-module $M$ and an $S$-module $N$, we say that $d(M,N)\leq 2^{-n}$, if $\app Rn$ and $\app Sn$ have a common scalar extension $T$ such that $M\tensor_RT\iso N\tensor_ST$. In particular, $d(R,S)\leq d(M,N)$. We will not study the similarity problem for modules---and at present, I do not know whether this is a smooth classification problem, even over a fixed ring. We will use this metric in the proof of Theorem~\ref{T:Poinc}; see also \cite[\S11]{SchFinEmb} for some further applications.

We now turn to some other classification problems that can be reduced to the classification up to similarity.

\subsection*{Classification of analytic germs}
Let $\kappa $ be a field. By an \emph{analytic germ} over $\kappa $, we mean a complete Noetherian local ring with residue field $\kappa $; we denote the set of isomorphism classes of analytic germs by $\siman(\kappa)$. Note that if $\kappa $ has prime \ch\ $p$, then the germ can either have equal  or mixed  \ch. By the Cohen structure theorem,  analytic germs are simply homomorphic images of power series rings $\pow VX$, with $V$ either $\kappa $ (equal \ch) or the complete $p$-ring over $\kappa $ (mixed \ch). Assume, moreover, that $\kappa $ is algebraically closed and has size of the continuum. It follows that every (countable) ultrapower $\ul \kappa $ of $\kappa $ is again algebraically closed and has the same cardinality as $\kappa $, whence by Leibniz's theorem, is  isomorphic with $\kappa $. In the mixed \ch\ case, by uniqueness of $p$-rings, the catapower of $V$ is then also isomorphic to $V$. This shows that   the set of analytic germs over such a field $\kappa $ is, up to isomorphism, closed under cataproducts, whence under limits. Moreover,  there are, up to isomorphism, only countably many analytic germs of dimension zero, and they form a dense subset $\siman_0(\kappa)$ of $\siman(\kappa)$. In conclusion, we showed Theorem~\ref{T:formgerm} from the introduction. Note that in the above, we may replace the size of the continuum by any cardinal of the form $2^\gamma$, with $\gamma$ an infinite cardinal. Moreover, under the Generalized Continuum Hypothesis, this means any uncountable cardinal.

\subsection*{Infinitesimal deformations}
By a \emph{\defpair} $\tuple R$, we mean a pair $(R,\tuple x)$, with $(R,\maxim)$ a Noetherian local ring and $\tuple x:=\rij xd$ a tuple generating an $\maxim$-primary ideal. To emphasize the maximal ideal, we may also represent the \defpair\ $\tuple R$ as the triple $(R,\maxim,\tuple x)$. We call $R$ and $\tuple x$  respectively the \emph{underlying ring} and \emph{tuple} of $\tuple R$, and we call the length of $R/\tuple xR$ the \emph{colength} of $\tuple R$. We call $\tuple R$ \emph{parametric}, if $\tuple x$ is a system of parameters. When  we say that a \defpair\ has a certain ring theoretic property, then we mean that its underlying ring has this property.  Let $\tuple S:=(S,\tuple y)$ be a second \defpair, with $\tuple y=\rij ye$. A \emph{morphism $\tuple R\to \tuple S$ of \defpair s},  is a ring \homo\ $R\to S$ mapping $\tuple x$ to $\tuple y$. In particular, there are no morphisms between \defpair{s} with tuples of different length. It is easy to verify  that these definitions make the class of \defpair{s} into a  category. We call a morphism $\tuple R\to \tuple S$ \emph{flat, unramified, a scalar extension, etc.}, \iff\ the underlying \homo\ $R\to S$ has this property. We say that $\tuple R$ and $\tuple S$ are \emph{similar}, in symbols, $\tuple R\approx\tuple S$, if they have a common scalar extension $\tuple T$ (as \defpair{s}). As before, we denote the similarity class of a \defpair\ $\tuple R$ by $[\tuple R]$.

The \emph{$n$-th infinitesimal deformation} of a \defpair\ $\tuple R:=(R,\tuple x)$, denoted $\app{\tuple R}n$, is by definition the Artinian \defpair\ $(R/\tuplepow xnR,\tuple x)$, where for an arbitrary tuple $\tuple y:=\rij ys$, we write  $\tuplepow yn$ for the tuple $(y_1^n,\dots,y_s^n)$. If $\tuple R\to \tuple S$ is a morphism of \defpair{s}, then it induces, for each $n$, a morphism $\app{\tuple R}n\to \app{\tuple S}n$.

\begin{lemma}\label{L:defpairscal}
If $\tuple R$ and $\tuple S$ are similar \defpair{}s, then $\app{\tuple R}n\approx\app{\tuple S}n$, for all $n$.
\end{lemma}
\begin{proof}
Since the respective underlying rings $R$ and $S$ are similar, they have the same dimension. Without loss of generality, we may assume that $\tuple R\to \tuple S$ is a scalar extension. By definition of morphism,  under the scalar extension $R\to S$, the tuple of $\tuple R$ is sent to that of $\tuple S$, and the assertion is now clear.
\end{proof}

Let $\simlocdef$ denote the set of similarity classes of \defpair{s} (the argument that this is indeed a set is analoguos to the one for $\simloc$). We define the \emph{deformation metric} on $\simlocdef$ in analogy with the \deform\ metric: given two similarity classes of \defpair{s} $[\tuple R:=(R,\tuple x)]$ and $[\tuple S:=(S,\tuple y)]$, we set $d([\tuple R],[\tuple S])\leq 2^{-n}$, if   $\app{\tuple R}n\approx\app{\tuple S}n$. By Lemmas~\ref{L:app} and \ref{L:defpairscal}, this definition is independent from the choice of representatives. Moreover, if $\app{\tuple R}n\approx\app{\tuple S}n$, then the definition of morphisms in the category of \defpair{s} implies that $\app{\tuple R}i\approx\app{\tuple S}i$, for all $i\leq n$. Indeed, we may reduce to the case that we have a scalar extension $\app{\tuple R}n\to\app{\tuple S}n$, which therefore maps $\tuple x$ to $\tuple y$, and the claim is now clear. If $d(\tuple R,\tuple S)<1$, then $\tuple R$ and $\tuple S$ have in particular the same colength.
As with rings, we will often indentify a similarity class with any \defpair\ contained in it, and so we will omit brackets in our notation and speak of the distance between \defpair{s}. The connection between the \deform\ metric and the deformational metric is given by:

\begin{proposition}\label{P:metcomp}
If  $\tuple R$ and $\tuple S$ are \defpair{s} with respective underlying rings $R$ and $S$, then  $d(R,S)\leq d(\tuple R,\tuple S)$. Conversely, for every \defpair\ $\tuple R$,   if $T$ is a Noetherian local ring at distance $\varepsilon$ from $R$, then we can find a \defpair\ $\tuple T$ with underlying ring $T$, such that $d(\tuple R,\tuple T)\leq\varepsilon^{1/{(lm+1)}}$, where   $l$ is the colength of $\tuple R$ and $m$ the length of its tuple. If, moreover, $\tuple R$ is parametric, and   $\op{dim}(R)=\op{dim}(T)$, then we may also choose $\tuple T$ to be parametric.
\end{proposition}
\begin{proof}
Let $(R,\maxim)$ and $(S,\mathfrak n)$ be the respective underlying rings of $\tuple R$ and $\tuple S$, and let $\tuple x$ and $\tuple y$ be their respective tuples. If $\app{\tuple R}k\approx\app{\tuple S}k$, for some $k$, then clearly $\app Rk\approx\app Sk$, since $\tuplepow xkR\sub\maxim^k$ and  $\tuplepow ykS\sub\mathfrak n^k$. This proves the first assertion. 

To prove the second, observe that since $\maxim^l\sub I:=\tuple xR$, we get 
\begin{equation}\label{eq:sublmn}
\maxim^{lmn}\sub I^{mn}\sub\tuplepow xn R,
\end{equation}
 for all $n$. Hence $\bar R:=R/\tuplepow xnR$ is a homomorphic image of $\app R{lmn}$. Suppose $d(R,T)\leq 2^{-k}$, so that $\app Rk\approx \app Tk$. Without loss of generality, we may assume that $\app Rk\to \app Tk$ is a scalar extension. Let $\tuple z$ be a lifting in $T$ of the image  of $\tuple x$  in $\app Tk$ under this scalar extension, and put $\tuple T:=(T,\tuple z)$. Let  $n$ be an integer strictly less than $k/lm$, so that $lmn<k$. We want to show that $\pr^k\sub \tuplepow znT$, where $\pr$ is the maximal ideal of $T$. Put $\bar T:= T/\tuplepow znT$. The map $\app Rk\to \app Tk$ induces a scalar extension $\bar R\to \bar T/\pr^k\bar T$. By \eqref{eq:sublmn},  the latter is annihilated by $\pr^{lmn}$.  Hence $\pr^{lmn}\bar T=\pr^k\bar T$, and since $lmn<k$, Nakayama's Lemma yields $\pr^{lmn}\bar T=0$, and the claim follows. In particular, base change induces a scalar extension $\bar R\to \bar T$, and hence a scalar extension $\app{\tuple R}n\to \app{\tuple T}n$ of \defpair{s}, showing that $d(\tuple R,\tuple T)\leq 2^{-n}$, as we wanted to show.
\end{proof}

\begin{theorem}\label{T:defsmooth}
Classification of \defpair{s} up to similarity is smooth, or, more precisely, $\simlocdef$ is a Polish space.
\end{theorem}
\begin{proof}
Let $\seq{\tuple R}w$ be a Cauchy sequence in $\simlocdef$, and let $\seq Rw$ be the corresponding sequence of underlying rings. By Proposition~\ref{P:metcomp}, this latter sequence is also Cauchy, whence has  a limit in $\simloc$ by Theorem~\ref{T:simcomp}. In fact, we may take the cataproduct $\ulsep R$ of the $\seq Rw$  as a representative of this limit. Since all tuples in $\seq{\tuple R}w$ must have the same length, their ultraproduct yields a finite tuple in $\tuple x$ in $\ulsep R$. Moreover, almost all $\seq{\tuple R}w$   have the same colength, which, by \los, is then also the length of $\ulsep R/\tuple x\ulsep R$. In particular, $(\ulsep R,\tuple x)$ is a \defpair. The second part of  Proposition~\ref{P:metcomp} shows that it is the limit of the $\seq{\tuple R}w$. This proves that $\simlocdef$ is complete. Remains to show that the subset $\simlocdef_0$ of Artinian \defpair{s} is countable and dense. However, we   argued in the proof of Theorem~\ref{T:smooth} that each similarity class of an Artinian local ring $R$ contains a representative $R_0$ with a finitely generated residue field. Given any (finite) tuple $\tuple x$, we may choose $R_0$ so that it also contains $\tuple x$. From this it is easy to see that $\simlocdef_0$ is countable, and density is also immediate.
\end{proof}

We denote the subset of similarity classes of parametric \defpair{s} by $\simlocpar$.  Dimension, as this is equal to the length of the tuple, partitions this space in the pieces $\simlocpar_d$. It follows immediately from the above proof   that each $\simlocpar_d$ is a complete subspace of $\simlocdef$. In particular, $\simloc_0$ is isometric with $\simlocpar_0$. However,   for $d>0$, it is no longer clear whether   $\simlocpar_d$ has a countable dense subset, and therefore, it might fail to be a Polish subspace.

\subsection*{Classification of polarized schemes up to isomorphism}
Our next application is to the classification of  
projective schemes. We will tacitly assume that a  \emph{projective
scheme} $X$
is always of finite type over some field $\kappa$.   A
\emph{polarization} of $X$ is a choice of a very ample line bundle $\mathcal L$ on
$X$; we refer to this situation also by calling $\mathfrak X:=(X,\mathcal L)$ a
\emph{polarized scheme over $\kappa$}, and  we say that $X$ is  the \emph{underlying projective scheme} of $\mathfrak X$. In particular, a polarization   $(X,\mathcal
L)$
corresponds to a closed immersion $i\colon X\to \mathbb P_\kappa^n$, for some $n$, 
such that
$\mathcal L\iso i^*\loc(1)$, where $\loc(1)$ is the canonical twisting sheaf
on $\mathbb P_\kappa^n$.

The
\emph{section} ring of a polarized
scheme $\mathfrak X:=(X,\mathcal L)$ is defined as the graded $\kappa$-algebra
$$
S(\mathfrak X):=\sum_{n=0}^\infty H^0(X,\mathcal L^n)
$$
Note that, since $\mathcal L$ is very ample, $S(\mathfrak X)$ is a \emph{standard graded $\kappa$-algebra}, meaning that it has no homogeneous components of negative degree, its degree zero component is $\kappa$, and, as an algebra over $\kappa$, it is generated by its homogeneous elements of degree one.    
 
The \emph{vertex algebra} of $\mathfrak X$ is the localization of $S(\mathfrak
X)$ at the irrelevant ideal of all elements of positive degree, and will be
denoted by $\op{Vert}(\mathfrak X)$. If $X$ is irreducible and reduced, then the field of
fractions of $S(\mathfrak X)$ (and hence of $\op{Vert}(\mathfrak X)$) is equal
to the function field $\kappa(X)$. In particular, $\op{Vert}(\mathfrak X)$ is
a birational invariant of $X$. In fact, more is true:
 the polarized scheme $\mathfrak X:=(X,\mathcal L)$
can be recovered from   its section ring $S:=S(\mathfrak X)$ as $X=\op{Proj}(S)$ and $\mathcal L=\widetilde{S(1)}$, where $S(1)$ is the \emph{Serre twist} of $S$. We therefore say that two polarized schemes $\mathfrak X:=(X,\mathcal L)$ and $\mathfrak Y:=(Y,\mathcal M)$ are \emph{isomorphic}, if their section rings are isomorphic as graded $\kappa$-algebras, and this is then equivalent with the existence of an isomorphism $f\colon X\to Y$ of projective schemes, such that $f^*\mathcal M=\mathcal L$.

Let
$\mathbb{Pol}_\kappa$ be the set of isomorphism classes of polarized
schemes over $\kappa$. We metrize this space via pull-back along the vertex functor, that is to say, 
$$
  d(\mathfrak X,\mathfrak Y):=  d(\op{Vert}(\mathfrak
X),\op{Vert}(\mathfrak Y)).
$$
The following
easy lemma allows us to calculate this distance function:

\begin{lemma}\label{L:polapp}
Let $\mathfrak X:=(X,\mathcal L)$ be a polarized scheme over $\kappa$ with   vertex
algebra $R:=\op{Vert}(\mathfrak X)$. For each
$n$, we have an isomorphism of graded Artinian $\kappa$-algebras
$$
\app Rn\iso \bigoplus_{i=0}^{n-1} H^0(X,\mathcal L^i).
$$
\end{lemma}
\begin{proof}
Let $S:=S(\mathfrak X)$ be the section ring of $\mathfrak X$, and let 
$\maxim$ be the irrelevant maximal ideal. Since $R=S_\maxim$, we have 
$\app Rn=S/\maxim^n$, for all $n$. Since $S$ is a standard graded algebra,
 $\maxim^n$ consists of all
elements of degree at least $n$, that is to say, 
 $$
\maxim^n=\bigoplus_{i\geq n} H^0(X,\mathcal L^i),
$$
from which the assertion follows immediately.
\end{proof}

We can now show that we have indeed a metric on $\mathbb{Pol}_\kappa$:

\begin{corollary}
If two polarized schemes $\mathfrak X$ and $\mathfrak Y$ over $\kappa$ are at distance zero, then they are isomorphic.
\end{corollary}
\begin{proof}
By Lemma~\ref{L:polapp}, their section rings are isomorphic, and we already argued that this means that the two polarized schemes are isomorphic.
\end{proof}

\subsubsection*{Proof of Theorem~\ref{T:pol}}
We will show that $\mathbb{Pol}_\kappa$ is a Polish space, and to this end, we need to show that it contains a countable dense subset and is closed under limits. By Lemma~\ref{L:polapp}, the polarizations  $(X,\mathcal L)$ of zero-dimensional projective schemes are dense. Any such scheme  is the base change of a zero-dimensional projective scheme $X_0$ over a finitely generated field, and since very ample line bundles are generated by their global sections, we may choose $X_0$ so that it admits a very ample line bundle $\mathcal L_0$ which induces the line bundle $\mathcal L$ on $X$ by base change. This shows that, up to isomorphism,  there are only countably many polarizations of zero-dimensional projective schemes.

So remains to show that every Cauchy sequence $\seq{\mathfrak X}w:=(\seq Xw,\seq{\mathcal L}w)$ in $\mathbb{Pol}_\kappa$ has a limit. Let $\seq Rw:=\op {Vert}(\mathfrak X)$, so that by definition, $\seq Rw$ is a Cauchy sequence of Noetherian local rings. Let $\ulsep R$ be the cataproduct of the $\seq Rw$. By the same argument as in the proof of Theorem~\ref{T:formgerm}, our assumption on the field $\kappa$ implies that $\ulsep R$ has residue field isomorphic to $\kappa$. Let $R$ be an isomorphic copy of $\ulsep R$ having residue field $\kappa$, and let $f\colon \ulsep R\to R$ be the corresponding isomorphism. Fix some $n$. By Lemma~\ref{L:polapp}, $\app Rn\iso \app{\seq Rw}n$, for all $w\gg 0$. In particular, the $n$-th homogeneous piece $S_n:=H^0(\seq Xw,\seq{\mathcal L}w^n)$ is independent from $w$, for $w$ sufficiently large. Since   $S_n$ has finite length, its ultrapower is equal to its catapower, and, therefore, via $f$, isomorphic to itself. Let $S:=\oplus_nS_n$.
One verifies that this is a standard graded $\kappa$-algebra. For instance, to define the ring structure on $S$, it suffices to define the multiplication of two homogeneous elements, say $a\in S_i$ and $b\in S_j$. Take $w$ large enough so that   $H^0(\seq Xw,\seq{\mathcal L}w^{i+j})$ is equal to $S_{i+j}$. Choose $\seq aw,\seq bw\in S(\seq{\mathfrak X}w)$ so that their images in $\seq Rw$ have ultraproducts $\ul a,\ul b\in \ul R$ with $f(\ul a)=a$ and $f(\ul b)=b$. By \los, $\seq aw$ and $\seq bw$ are homogeneous of degree $i$ and $j$ respectively. We then define $ab$ as the image under $f$ of the ultraproduct  of the $\seq  aw\seq b w\in S_{i+j}$. The other properties are checked similarly. In particular, $\app Rn\iso S_0\oplus\dots\oplus S_n$, showing that $R$ is the localization of $S$ at its irrelevant maximal ideal. Let $\mathfrak X:=(X,\mathcal L)$ be the   polarized scheme determined by $S$, namely, $X:=\op{Proj}(S)$ and $\mathcal L:=\widetilde{S(1)}$. Remains to show that $\mathfrak X$ is the limit of the $\seq{\mathfrak X}w$, and this is immediate form the fact that $\op{Vert}(\mathfrak X)=R$.\qed

\section{Prolegomena to a complete set of invariants: slopes}

Theorem~\ref{T:smooth}, although promising, is far from an efficient
classification up to similarity. In this final section, we will discuss some
(albeit feable) attempts to make it more concrete. As mentioned in the
introduction,
any two (uncountable) Polish spaces are Borel equivalent, namely to the standard
Borel space
on the reals. So, given any (concrete) Polish space $\mathbb B$, we ask for a
Borel bijection
$q\colon \simloc\to\mathbb B$. 

Let us call a map $q\colon \simloc\to\mathbb B$ a \emph{slope}, if it is
continuous. Of course, the identity map into $\simloc$ itself is a slope, but
we seek more concrete examples. A solution to the
classification
problem would, for instance,  be provided by any real-valued, injective slope. Extending this  terminology, let
us
 say that for some subset $\mathbb S\sub \simloc$ and  a map $q\colon  \simloc\to \mathbb B$, that $q$ is a \emph{slope}   on $\mathbb S$  when its restriction to
$\mathbb S$ is continuous. A priori, the
theory only predicts that we can find a real-valued, injective Borel map, which
in general is only continuous outside a meagre subset, but perhaps we may venture to postulate the existence of a countable partition $\{\mathbb S_i\}$ of $\simloc$, and an injective map $q\colon\simloc\to \mathbb R$, such that $q$ restricted to each piece $\mathbb S_i$ is a slope. Moreover, we want this partition to be indexed by some natural  discrete invariants that are preserved under the similarity relation, like dimension and/or parameter degree.    We start with some  examples of   non-injective slopes taking values into a concrete
complete topological space (from now on, we will confuse a similarity class
with any of its members):

\begin{proposition}\label{P:Hilb}
Viewing $\pow\zet t$ in its $t$-adic metric, the map $\op{Hilb}\colon\simloc\to
\pow\zet t$ induced by associating
to a Noetherian local ring its Hilbert series $\op{Hilb}(R)$, is a slope. 
\end{proposition}
\begin{proof}
Recall that the \emph{Hilbert series} of $(R,\maxim)$ is defined to be the formal
power series 
$$
\op{Hilb}(R):=\sum_{n=0}^\infty \ell(\maxim^n/\maxim^{n+1})t^n.
$$
If $R$ and $S$ are similar, then they have the same
Hilbert
series, showing that $\op{Hilb}$  is defined on $\simloc$. By an easy calculation,  
$\ell(\maxim^n/\maxim^{n+1})=\ell(\app
R{n+1})-\ell(\app Rn)$. Hence, if $\seq Rw$ converges to $R$, then for each $n$, we
have $\app
Rn=\app{\seq Rw}n$, for all sufficiently large $w$, showing that $\op{Hilb}$ is
continuous. 
\end{proof}

For a second example of a (non-injective) slope, we make the following definition. Let $R$ be a   local ring with residue field $\kappa $, and let $M$ be a finitely generated $R$-module. The \emph{$n$-th Betti number} $\beta_n(M)$ of $M$ is defined as the vector space dimension of $\tor RnM \kappa $. Alternatively, at least in the Noetherian case, the Betti numbers are the ranks in a minimal free resolution of $M$, and, hence by Nakayama's Lemma, the minimal number of generators of the syzygies of $M$. 
The generating series of these Betti numbers, that is to say, the formal power series
$$
\poinc M{}:=\sum_n\beta_n(M)t^n
$$
is called the \emph{Poincare series} of $M$.  We define the \emph{residual Poincare series} of $R$ to be the Poincare series of its residue field, and denote it $\poinc R{res}$. If $R\to S$ is a scalar extension, and $F_\bullet$ a minimal free resolution of the residue field $\kappa$ of $R$, then by flatness, $F_\bullet\tensor_RS$ is a minimal free resolution of $\kappa\tensor_RS$, and the latter is the residue field of $S$, since $R\to S$ is unramified. Hence, any two similar rings have the same residual Poincare series. Let $\simcm$   be the subset of $\simloc$ consisting of all similarity classes of   local \CM\ rings. By \cite[Corollary 8.7]{SchFinEmb},  if we also fix dimension and multiplicity, then each $\simcm_{d,e}$ is closed under limits.

\begin{theorem}\label{T:Poinc}
The residual  Poincare series is a slope on each class $\simcm_{d,e}$.
\end{theorem}
\begin{proof}
The continuity of the map associating to a $d$-dimensional local \CM\ ring $R$ of multiplicity $e$ its residual Poincare series $\poinc R{res}$ is an immediate consequence of \cite[Theorem 11.4]{SchFinEmb}. Indeed, if $\seq Rw$ is a Cauchy sequence, then, for any fixed $n$, the residue field of each $\seq R w$ has the same   $n$-th  Betti number, for $w$ sufficiently large, by the cited result. By \cite[Proposition 8.9]{SchFinEmb}, this is then also the Betti number of  the cataproduct $\ulsep R$, that is to say, up to similarity, the limit of the $\seq Rw$.  Therefore, $\poinc{\seq Rw}{res}$ converges, in the $t$-adic topology, to $\poinc{\ulsep R}{res}$.
\end{proof}

For a local  \CM\ ring $R$,   we also define its \emph{canonical Poincare series}, denoted $\poinc R{can}$, as the Poincare series $\poinc {\omega_{\complet R}}{}$ of the canonical module $\omega_{\complet R}$ of its completion $\complet R$ (note that the canonical module always exists when the ring is complete; see for instance \cite[\S3.3]{BH}).
In particular,  $R$ is Gorenstein \iff\ its canonical Poincare series is constant (equal to $1$): indeed, $R$ is Gorenstein \iff\ $\omega_{\complet R}\iso \complet R$.  It is not hard to check that the canonical Poincare series is independent from the choice of representative of a similarity class of a \CM\ local ring (by the same argument as in \cite[Theorem 3.3.5(c)]{BH}). 
Let $R$ and $S$ be two rings in $\simcm_{d,e}$ at distance at most $2^{-de-1}$. After a scalar extension, we may assume that $R$ and $S$ are complete, with infinite residue fields $\kappa $ and $\lambda$, respectively. In particular, there exists a system of parameters $\tuple x$ in $R$ such that $\bar R:=R/\tuple xR$ has length $e$.   Proposition~\ref{P:metcomp} then yields a system of parameters $\tuple y$ in $S$ such that $\bar R\approx \bar S:=S/\tuple yS$. Since the canonical module $\omega_R$ is maximal \CM,    $\tuple x$ is   $\omega_R$-regular, and likewise, $\tuple y$ is  $\omega_S$-regular. Therefore, 
\begin{equation}\label{eq:torcan}
\begin{aligned}
\tor Rn{\omega_R}\kappa &\iso \tor {\bar R}n{\omega_R/\tuple x\omega_R}\kappa\\
\tor Sn{\omega_S}\lambda &\iso \tor {\bar S}n{\omega_S/\tuple y\omega_S}\lambda
\end{aligned}
\end{equation}
for all $n$. Since $\bar R$ and $\bar S$ are similar, they have the same canonical Poincare series $P(t)$. By \cite[Theorem 3.3.5]{BH}, the respective canonical modules of $\bar R$ and $\bar S$ are $\omega_R/\tuple x\omega_R$ and $\omega_S/\tuple y\omega_S$. By \eqref{eq:torcan}, therefore, the canonical Poincare series of $R$ and $S$ are both equal to $P(t)$. In conclusion, we showed:

\begin{proposition}
On each ball of radius $2^{-de-1}$ in $\simcm_{d,e}$, the canonical Poincare series is constant. In particular, if one of its members is Gorenstein, then so is any.\qed
\end{proposition}

\subsection*{Quasi-slopes}
To find a slope, it is enough to have it defined on the countable dense open
subset $\simloc_0$ given by Theorem~\ref{T:smooth}.  For
any map $q_0\colon \simloc_0\to \mathbb B$ into a complete metric space (not
necessarily continuous),
define
its extension $\complet q_0$ as the partial map $\simloc\dasharrow  \mathbb B$
given as the limit
 of  the $q_0(\app Rn)$, for $R$ a Noetherian local ring, whenever this limit
exists. Note that if $R\approx S$, then $q_0(\app Rn)$ converges \iff\
$q_0(\app Sn)$ does, and their limits are similar. In particular, $q_0(R)=\complet q_0(R)$ whenever $R$ is
Artinian. We call a map $q_0\colon \simloc_0\to \mathbb B$ a
\emph{quasi-slope}, if $\complet q_0$ is everywhere defined. By abuse of
terminology, we then also refer to this extension $q:=\complet q_0$ as a
quasi-slope. In other words, $q\colon \simloc\to
\mathbb B$ is a quasi-slope
if, $q(\app Rn)$ converges to $q(R)$, for every Noetherian local ring $R$. 
The following
corollary is now immediate from Theorem~\ref{T:smooth}:

\begin{corollary}\label{C:slope}
Any continuous map $q_0\colon \simloc_0\to \mathbb B$ is a quasi-slope and its 
extension $\complet q_0$ is a
slope.\qed
\end{corollary}

We next   show how some of the usual invariants, although in general not
slopes, become quasi-slopes when properly modified.
Let $\delta_0$ be defined on $\simloc_0$ as follows. Given an
Artinian local ring $(A,\maxim)$, let $n$ be its degree of nilpotency (that is to
say, the least $k$ such that $\maxim^k=0$). Put
$$
\delta_0(A):=\log_2(\frac{\ell (A)}{\ell (\app A{n/2})})
$$
where for a postive real number $r$, we define $\app Ar:=\app Az$ with
$z:=\op{int}(r)$ 
the largest integer less than or equal to
$r$.  

\begin{proposition}
The map $\delta_0$ is a quasi-slope. In fact,  $\complet\delta_0(R)$ is equal to
the dimension of $R$, whenever this dimension is non-zero.
\end{proposition}
\begin{proof}
Let $R$ be a Noetherian local ring of dimension $d>0$. By the Hilbert-Samuel
theory, there exists a polynomial $P_R\in\pol{\mathbb Q}t$
of degree $d$, such that $\ell(\app Rn)=P_R(n)$ for $n\gg 0$. Hence $\delta_0(\app
Rn)=\log_2(P_R(n)/P_R(\op{int}(n/2)))$, for $n\gg 0$. It is now an exercise to
show that for any polynomial $P$ of degree $d$, the limit of $P(n)/P(\op{int}(n/2))$ is equal to $2^d$.
\end{proof}

In view of this result, we call $\complet\delta_0(R)$ the
\emph{quasi-dimension} of $R$. So only Artinian local rings have a
quasi-dimension which is different from their (Krull) dimension.
Using the formula
$$
\lim_{n\to \infty} \frac{P(\op{int}(\sqrt n))^2}{P(n)}=a_d
$$
where $a_d$ is the leading coefficient of a polynomial $P$, we get by a
similar argument that the map $\epsilon_0$ defined on $\simloc_0$ by the
condition $\epsilon_0(A):=(\ell(\app A{\sqrt n}))^2/\ell(A)$ is a quasi-slope,
and $\complet\epsilon_0 (R)=e/d!$ whenever $d>0$, where
$e$ is the multiplicity of $R$ and $d$ its dimension.

Several questions now arise naturally: what is the nature of the subset of
$\simloc$ of all Noetherian local rings of a fixed quasi-slope? Can we break up
(or even stratify) $\simloc$ in ``natural''
pieces on which a quasi-slope becomes  continuous. I will conclude with an
example
of how one can answer the second question for quasi-dimension.  For
a Noetherian local ring $R$, define $\rho(R)$ as the
supremum over all $n$ of
$$
\norm{\frac{d!\ell(\app Rn)}{en^{d-1}}-n}
$$
where $d:=\op{dim}(R)$ and $e:=\op{mult}(R)$ are respectively the dimension and
multiplicity of $R$.
In other words, $\rho(R)$ is the smallest real number $\rho\geq 0$ such that
\begin{equation}\label{eq:rho}
n^d-\rho n^{d-1} \leq \frac{d!}e\ell(\app Rn)\leq  n^d+\rho
n^{d-1}
\end{equation}
for all $n>0$. That this supremum exists is an easy consequence of the
Hilbert-Samuel theory. For instance, if $R$ is Artinian of length $l$, then
$1\leq \rho(R)\leq
l$, but these bounds are not sharp. Note that $\rho$ is \emph{not} a quasi-slope
(this is easily checked for $R:=\pow \kappa t$).

This new invariant determines the rate of
convergence in the definition of the quasi-dimension, as the next result shows.
To formulate it, we use $\near r$ to denote the rounding to the nearest integer
of a real number $r$, that is to say, $\near r$ is the unique integer inside the
half open interval $[r-\frac 12,r+\frac 12)$.

\begin{lemma}
For a  Noetherian local ring $R$, if $n\geq 10\rho(R)$, then
$\near{\delta_0(\app Rn)}$ is equal to its dimension.
\end{lemma}
\begin{proof}
Let $b:=\rho(R)$, and let $d$ and $e$ be the respective dimension and
multiplicity of $R$.
Using
\eqref{eq:rho}, we get estimates 
\begin{equation}\label{eq:bn}
1-\frac bn\leq\frac {d!\ell(\app Rn)}{en^d}\leq 1+\frac bn
\end{equation}
for all $n$. In the convergence of $\delta_0$ we may take the limit over even
$n$ only, so let us assume that $n=2m$. Dividing inequalities \eqref{eq:bn} for $2m$ by those for $m$,   we get estimates
$$
 2^d\big(\frac {1-\frac b{2m}} {1+\frac bm}\big)\leq\frac
{\ell(\app R{2m})}{\ell(\app Rm)}= 2^{\delta_0(\app Rn)} \leq 2^d\big(\frac {1+\frac b{2m}} {1-\frac
bm}\big).
$$
Hence, if the ratio between the two outside fractions is strictly less than
$2$, then after taking the logarithm with base two, they become the endpoints
of an interval $[\alpha,\beta]$ 
of length stricly less than one, containing $\delta_0(\app Rn)$. Since $\alpha<d<\beta$, the only integer in $[\alpha,\beta]$ is  $d$, showing that $\near{\delta_0(\app Rn)}=d$.

 For the ratio to be at most $2$, we need
$$
(m+\frac b2)(m+b)<2(m-\frac b2)(m-b)
$$
and a simple calculation shows that this is true whenever $m>5b$.
\end{proof}

Immediately from this we get:

\begin{corollary}
For each $b\in\nat$, let $\simloc_{\rho\leq b}$ be the subset of $\simloc$
consisting of all Noetherian local rings $R$ such that $\rho(R)\leq b$, thus
yielding a
filtration $\simloc_{\rho\leq 0}\sub \simloc_{\rho\leq 1}\sub\dots$ of
$\simloc$.
Then quasi-dimension is a slope on each $\simloc_{\rho\leq b}$. 
\end{corollary}

A similar argument can be used to show that $\epsilon_0$ is a slope on each
$\simloc_{\rho\leq b}$, by establishing an analogous bound for the convergence of
$\epsilon_0$ to $e/d!$ which only depends on $d:=\op{dim}(R)$, $e:=\op{mult}(R)$
and $\rho=:\rho(R)$. As far
as bounding $\rho$ itself  is concerned, if $R$ is \CM, then it is bounded as
a function
of $e$ and $d$ only, but without the \CM\ assumption this is probably false. In
the latter case, we can use any ``big degree'' $D$ \`a la Vasconcelos to arrive
at such a bound in terms of $d$ and $D(R)$ (this is an easy consequence of  
\cite[Theorem 4.1]{RTV}).

Recall (see, for instance, \cite[III. Ex. 5.1]{Hart}) that
the
\emph{Euler characteristic} of a projective scheme $X$  is
defined by the formula
$$
\chi(X):=\sum_i (-1)^ih^i(X,\loc_X)
$$
where $h^i(X,\loc_X)$ is the vector space dimension of the sheaf cohomology
$H^i(X,\loc_X)$. In particular, if $X$ is a curve, then $\chi(X)-1$ is the
\emph{genus} of $X$.

\begin{proposition}
The map $\mathbb{Pol}\to \mathbb R$ sending a polarized scheme $\mathfrak
X=(X,\mathcal L)$ to $\chi(X)$ is continuous. 
\end{proposition}
\begin{proof}
We calculate the Euler characteristic by means of the Hilbert-Samuel polynomial
$P_X(n)$ of
$X$ as 
\begin{equation}\label{eq:EC}
\chi(X)=P_X(0).
\end{equation}
 In fact, as it is a birational invariant, we may calculate the Euler characteristic by means of any
polarization $\mathfrak
X=(X,\mathcal L)$  of $X$. The \emph{Hilbert series} of $\mathfrak X$ is
defined as
$$
\op{Hilb}(X):=\sum_{n=0}^\infty h^0(X,\mathcal L^n)t^n.
$$
Let $(R,\maxim):=\op{Vert}(\mathfrak X)$. It is not hard to see that
$H^0(X,\mathcal L^n)\iso \maxim^n/\maxim^{n+1}$ and hence that
$X$ and $R$ have the same Hilbert series and the same Hilbert-Samuel polynomial. Moreover, the connection between the Hilbert series $h(t)$
  and the corresponding Hilbert-Samuel polynomial $P(n)$ is
given by the formula
\begin{equation}\label{eq:hilb}
P(n)= \sum_{j=0}^{d-1} \frac{(-1)^j}{j!} \binomial{n+d-1-j}n\restrict{\frac{\partial^j}{\partial
t}\big((1-t)^dh(t)\big)}{t=1},
\end{equation}
where $d$ is the degree of $P$ (that is to say, the dimension of $X$).
Moreover, if $h_i$ are Hilbert series with corresponding Hilbert polynomial
$P_i$, then from
the fact that $$
\restrict{\frac{\partial^j}{\partial t}\big((1-t)^dt^n\big)}{t=1} =0
$$
for all $j<d$ and all $n$, we get from \eqref{eq:hilb} that $P_1=P_2$ whenever
$h_1$ and $h_2$ are $t$-adically close.
Therefore, by Proposition~\ref{P:Hilb} and \eqref{eq:EC}, the Euler
characteristic is continuous.
\end{proof}

\subsection*{Future work}
In work in progress (\cite{SchSchGr}), we assign to any Artinian  $\kappa$-algebra, with $\kappa$ an \acf\ of size the continuum, a first-order formula modulo the theory of Artinian $\kappa$-algebras. Associating to this  theory   its Grothendieck ring $K_0:=K_0(\kappa)$ on pp-formulae, we get a map $\siman_0(\kappa)\to K_0\colon R\mapsto [ R]$, which is compatible with direct sum and tensor product, and which becomes injective when we replace the isomorphism relation with a stable version of it. Hence, we may associate to any analytic germ $R$, its \emph{formal Hilbert series}
$$
\op{Hilb}^{\op{form}}(R) :=\sum_n  [\app Rn] t^n\in\pow {K_0}t
$$
This would yield a complete invariant modulo the Grothendieck ring $K_0$ (this Grothendieck ring, however, admits the classical Grothendieck ring of $\kappa$ as a homomorphic image, whence is potentially a very complicated object).

 \providecommand{\bysame}{\leavevmode\hbox to3em{\hrulefill}\thinspace}
\providecommand{\MR}{\relax\ifhmode\unskip\space\fi MR }
\providecommand{\MRhref}[2]{%
  \href{http://www.ams.org/mathscinet-getitem?mr=#1}{#2}
}
\providecommand{\href}[2]{#2}

\end{document}